\documentclass[11pt,a4]{article}

\usepackage{ae} 
\usepackage[T1]{fontenc}
\usepackage[ansinew]{inputenc}
\usepackage{amsmath,amssymb,amsthm}
\usepackage{graphicx}
\usepackage{bbding}
\usepackage{color}
\usepackage[colorlinks,citecolor=red, linkcolor=blue,hyperindex,CJKbookmarks]{hyperref}
\usepackage{lscape}
\usepackage{epsfig}
\usepackage{mathrsfs}
\usepackage{dsfont}
\usepackage{bbm}
\usepackage{fancyhdr}
\usepackage{multirow}

\hypersetup{pdfpagemode=FullScreen}

\newtheorem{theorem}{Theorem}[section]
\newtheorem{corollary}{Corollary}[section]

\newtheorem{lemma}{Lemma}[section]

\theoremstyle{definition}
\newtheorem{definition}[theorem]{Definition}
\newtheorem{remark}{Remark}[section]

\newtheorem{example}{\hskip\parindent\bf Example}[section]
\numberwithin{equation}{section}

\newcommand{\al}{\alpha}
\newcommand{\bt}{\beta}
\newcommand{\ba}{\begin{array} }
\newcommand{\ea}{\end{array} }
\newcommand{\be}{\begin{equation} }
\newcommand{\ee}{\end{equation} }
\newcommand{\baa}{\begin{align} }
\newcommand{\eaa}{\end{align} }
\newcommand{\da}{\delta}
\newcommand{\kal}{\kappa}

\newcommand{\la}{\lambda}
\newcommand{\f}{\displaystyle\frac}
\newcommand{\il}{\displaystyle\int}

\newcommand{\vp}{\varphi}

\newcommand{\ga}{\gamma}
\newcommand{\ep}{\epsilon}

\newcommand{\ra}{\rightarrow}
\newcommand{\Og}{\Omega}
\newcommand{\og}{\omega}

\newcommand{\htt}{\widehat{T}}
\newcommand{\wpt}{\widehat{P}}
\newcommand{\wqt}{\widehat{Q}}
\newcommand{\wkt}{\widehat{K}}
\newcommand{\pt}{\widetilde{P}}
\newcommand{\qtt}{\widetilde{Q}}
\newcommand{\wa}{\widetilde{A}}
\newcommand{\wb}{\widetilde{B}}

\newcommand{\sg}{\sigma}

\newcommand{\R}{{\mathbb{R}}}

\newcommand{\W}{{\mathcal{W}}}
\newcommand{\Xa}{{\mathcal{X}}}

\newcommand{\B}{{\mathcal{B}}}

\newcommand{\tx}[1]{\quad\mbox{#1}\quad}

\newcommand{\tu}[1]{\textup{#1}}

\newcommand{\ve}{\varepsilon}
\newcommand{\iy}{\infty}

\newcommand{\lf}{\left(}
\newcommand{\rf}{\right)}
\newcommand{\hc}{\hat{c}}
\DeclareMathOperator{\id}{id} \DeclareMathOperator{\Imm}{Im}

\begin{document}
\date{}
\title{Nonuniform $(h,k,\mu,\nu)$-Dichotomy with Applications to
Nonautonomous Dynamical Systems
\thanks{MF is supported by NSFC-11271065 and RFPD-20130043110001; JZ and YL are
supported by NSFC-11201128, NSFHLJ-A201414, CPSF,
STIT-HEI-HLJ-2014TD005, HLJUF-DYS-JCL201203.}}
\author{Jimin Zhang$^a$, Meng Fan$^b$\thanks{Corresponding
author. Email: mfan@nenu.edu.cn}, Liu Yang$^a$} \maketitle

\begin{center}
\begin{minipage}{12cm}
\begin{description}

\item \small $a$.~School of Mathematical Sciences, Heilongjiang
University, 74 Xuefu Street, Harbin, Heilongjiang, 150080, P. R.
China

\item \small $b$.~School of Mathematics and Statistics, Northeast
Normal University, 5268 Renmin Street, Changchun, Jilin, 130024,
P. R. China

\end{description}
\end{minipage}
\end{center}

\vskip 0.2in

\begin{abstract}
The paper develops and studies a very general notion of dichotomy, referred to as
``nonuniform $(h,k,\mu,\nu)$-dichotomy". The new notion contains as special cases
most versions of dichotomy existing in the literature.
The paper then provides corresponding new versions of robustness,
Hartman-Grobman theorem, and stable manifold theorem for nonautonomous dynamical systems
in Banach spaces in term of the nonuniform $(h,k,\mu,\nu)$-dichotomy.
\end{abstract}

\section{Introduction}

The well-known and established notion of dichotomy in the linear
analysis of nonautonomous systems, essentially originated in
landmark work of Perron \cite{Perron1930}, extends the idea of
hyperbolicity from autonomous systems to explicitly time-dependent
ones. Dichotomy is of fundamental importance and plays a central
role in nonautonomous dynamics. It entails a clear and simple
structure of the (extended) phase space, which in turn can be
utilized to address important questions regarding, for instance,
robustness, topological conjugacy, and invariant manifolds of nonautnomous dynamical systems (see Table
\ref{tabletypedichotomy} for the a summary of recent progress in the theory and applications of dichotomies). In the classical words of Coppel
\cite{Coppelbook1978}, ``dichotomies, rather than Lyapunov's
characteristic exponents, are the key to questions of asymptotic
behavior for nonautonomous differential equations". Dichotomies
have been the subject of extensive research over (at least) the
past three decades, leading to exciting new results in areas as
diverse as skew product flows, nonuniformly hyperbolic systems,
and finite-time dynamics. We refer the reader to the references
cited in Table \ref{tabletypedichotomy} for more details.

\begin{table}[h]
\begin{center}
\caption[]{Type of dichotomies and
their applications}\label{tabletypedichotomy} \scriptsize
\begin{tabular}{|p{1.5cm}|p{1.6cm}|p{2cm}|p{2cm}|p{2cm}|p{1.6cm}|}\hline
{\bf Dichotomy} &{\bf Definition}&{\bf Existence}&{\bf
Robustnesss} &{\bf Topological conjugacy} &{\bf Invariant
manifolds}\\\hline {\bf
UED}&{\bf\cite{Coppelbook1978,Megan2003,Perron1930}}&
{\bf\cite{Chow1995,Coppelbook1978,Latushkin1998,Perron1930,Preda2006,Ruan2005,Sacker1974}}
&{\bf\cite{Chow1995,Ju2001,Mendez2008,nrmp1998,Pliss1999,Popescu2009,plh2006}}
&{\bf\cite{Cong2007,Palmer1979,Palmer1973,Popescu2009,Popescu2004,Pugh1969,Shirikyan2004}}
&{\bf\cite{Chow1991JDE,Duan2003,Shirikyan2004}}\\\hline {\bf
$(h,k)$-D} &{\bf\cite{Naulin1997,Naulin1995,Pinto1992}}&
{\bf\cite{Naulin1995}}&{\bf \cite{Naulin1997,Naulin1995}} &{\bf}
&{\bf\cite{Fenner1997}}\\\hline {\bf NUED} &{\bf
\cite{Barreira2008book}}
&{\bf\cite{Barreira2008book}}&{\bf\cite{Barreira200855}}
&{\bf\cite{Barreira2006}} &{\bf\cite{Barreira2005d}}\\\hline {\bf
NUPD} &{\bf\cite{Barreira2009a}}
&{\bf\cite{Barreira2009a}}&{\bf\cite{bfvz2011}} &{\bf
}&{\bf\cite{bfvz2010,bfvz2011b}}\\\hline {\bf
N-$(\mu,\nu)$-D}&{\bf \cite{Bento2013,Chang2011}
}&{\bf\cite{Barreira2013}}&{\bf \cite{Chang2011,Chu2013}}&{}
&{\bf\cite{Bento2012,Bento2013}}\\
\hline {\bf $\rho$-NUED}&{\bf \cite{Barreira2008d,Barreira2009b}
}&\bf\cite{Barreira2008d}&{\bf \cite{Barreira2009b}
}&{\bf\cite{Barreira2009d} }&{\bf\cite{Barreira2009e}}\\\hline
\end{tabular}
$$\begin{array}{rll}
\hbox{Here, }&\hbox{UED: Uniform exponential
dichotomy}&\hbox{$(h,k)$-D:
$(h,k)$ dichotomy}\\
&\hbox{NUED: Nonuniform exponential dichotomy}&\hbox{NUPD:
Nonuniform polynomial dichotomy}\\&\hbox{N-$(\mu,\nu)$-D:
Nonuniform
$(\mu,\nu)$-dichotomy}&\hbox{$\rho$-NUED:$\rho$-nonuniform
exponential dichotomy}
\end{array}
$$
\end{center}
\end{table}

However, dynamical systems exhibit various different kinds of
dichotomic behavior and the existing notions of dichotomic type are too restrictive and can not well relate all those dichotomic behavior. It is important and of great
interest to look for more general hyperbolic behavior
\cite{Barreira2008book,Barreira2009b,Barreira2009a,Bento2013,Bento2014,Chang2011,Chu2013,Fenner1997,Megan2002,Minda2011,Naulin1995}.

In this paper, we
introduce a new notion called nonuniform
$(h,k,\mu,\nu)$-dichotomy,  which is a more general framework and includes as particular cases most versions of uniform and nonuniform dichotomy
usually found in the literature. The new notion's wide scope is
due to the fact that it allows different growth rates in the
stable space and unstable space and in the uniform part and
nonuniform part, and the comparison functions $h, k, \mu, \nu$ are
only assumed to be increasing and unbounded, but no specific form
(e.g. exponential or polynomial) is prescribed for them.

As revealed in Section \ref{existencesection}, where the existence
of nonuniform $(h,k,\mu,\nu)$-dichotomy is characterized in terms
of Lyapunov exponents and Lyapunov functions for a linear nonautnonomous dynamical system
in a finite-dimensional space, the notion of nonuniform
$(h,k,\mu,\nu)$-dichotomy is not just a routine extension of the
extensive notions of uniform or nonuniform dichotomies, but arises
naturally and has a more comprehensive description for the
qualitative and stable behavior of linear nonautonomous dynamical
systems.

The principal aim of the present paper, then, is to provide
corresponding new versions of the robustness (Section
\ref{robustnesssection}), Hartman-Grobman theorem (Section
\ref{ghsection}), and stable manifold theorem (Section
\ref{manifoldsection}) for nonautonomous continuous dynamical systems in Banach space.
The study reveals that  the new defined dichotomy still allows us to obtain results that generalize the ones in the literature.

\section{Definition of nonuniform $(h,k,\mu,\nu)$-dichotomy
}\label{definitionsection} \noindent

Let $\B(X)$ be the space of bounded linear operators on a Banach
space $X$. Consider the linear system
\begin{equation}\label{eqlinear}
x'=A(t)x,~~t\in\R,
\end{equation}
where $A(t)\in \B(X)$. Let $T(t,s)$ be the evolution
operator of \eqref{eqlinear} satisfying $T(t,s)x(s)=x(t), t, s\in\R$ for any
solution $x(t)$ of \eqref{eqlinear}. An increasing function $u:
\R\to (0,+\infty)$ is said to be a growth rate if $u(0)=1$,
$\lim\limits_{t\ra\iy}u(t)=\iy$ and $\lim\limits_{t\ra-\iy}u(t)=0$.
In the following, we always assume that $h(t), k(t), \mu(t), \nu(t)$ are growth
rates.

\begin{definition}\label{dehk}
\eqref{eqlinear} is said to admit a \emph{nonuniform
$(h,k,\mu,\nu)$-dichotomy} on $\R$ if there exists a projection
$P(t)$ such that
$
P(t)T(t,s)=T(t,s)P(s), \quad t, s \in\R
$
and there exist constants $a<0\leq b$, $\ve\ge0$ and $K>0$ such
that
\begin{equation}\label{deeqaaa}
\begin{split}
\|T(t,s)P(s)\|& \le K\left(\f{h(t)}{h(s)}\right)^a\mu(|s|)^\ve,~t\geq s,\\
\|T(t,s)Q(s)\|& \le
K\left(\f{k(s)}{k(t)}\right)^{-b}\nu(|s|)^\ve,~s\geq t,\end{split}
\end{equation}
where $Q(t)=\id-P(t)$ are the complementary projections of $P(t)$.
\end{definition}

In Definition \ref{dehk}, four different functions or growth rates
$h,k,\mu,\nu$ are chosen to characterize the stable space and the
unstable space and also the uniform part and the nonuniform part.
Intuitively, on the level of generality adopted here, the
parameters $a, b$ and $\ve$ can be made part of $h, k, \mu$ and
$\nu$, at least in the non-trivial case $a\ve<0$ and this would
simplify many expressions later on. For example, when $a\ve<0$,
one can replace \eqref{deeqaaa} with
$$
\ba{l}
\|T(t,s)P(s)\|\le K\left(\f{\bar{h}(t)}{\bar{h}(s)}\right)^{-1}\bar{\mu}(|s|),~t\geq s,\\
\|T(t,s)Q(s)\|\le
K\left(\f{\bar{k}(s)}{\bar{k}(t)}\right)^{-1}\bar{\nu}(|s|),~s\geq
t, \ea
$$
where $\bar{h}(t), \bar{k}(t), \bar{\mu}(t), \bar{\nu}(t)$ are
growth rates, or more general functions $a(t,s), b(t,s)$ (see \cite{Bento2014}, Bento and Silva consider two general  bounded functions on $\R^+$). In the present paper, we deliberately prefer to use
\eqref{deeqaaa}. The reason is that $a$ and $b$ play the role of
Lyapunov exponents while $\ve$ measures the nonuniformity of
dichotomies and the nonuniform
$(h,k,\mu,\nu)$-dichotomy can be more closely connected with Lyapunov exponents,  Lyapunov functions, and the
theory of nonuniform hyperbolicity (see Section
\ref{existencesection}). This also implies that the nonuniform $(h,k,\mu,\nu)$-dichotomy
exists widely in the theory of linear nonautonomous dynamical systems.

The nonuniform $(h,k,\mu,\nu)$-dichotomy is much more general and
extends the existing notions of uniform and nonuniform dichotomies such as uniform exponential dichotomy
($h(t)=k(t)=e^t, \ve=0$) \cite{Coppelbook1978}, $(h,h)$-dichotomy
($h(t)=k(t), \ve=0$) \cite{Pinto1992},  $(h,k)$-dichotomy
($\ve=0$) \cite{Naulin1995}, nonuniform exponential dichotomy
($h(t)=k(t)=e^t,\mu(t)=\nu(t)=e^{|t|}$) \cite{Barreira2008book,Zhang2014},
nonuniform polynomial dichotomy ($h(t)=k(t)=\mu(t)=\nu(t)=t+1$ for $t\in\R^+$)
\cite{Barreira2009a}, nonuniform $(\mu,\nu)$-dichotomy
($h(t)=k(t)=\mu(t)$ and $\mu(t)=\nu(t)=\nu(t)$ for $t\in\R^+$) \cite{Chang2011,Bento2013,Barreira2013} ,
$\rho$-nonuniform exponential dichotomy
($h(t)=k(t)=\mu(t)=\nu(t)=e^{\rho(t)}$ for $t\in\R^+$) \cite{Barreira2009b} and
so on.

The following contrived example shows the generality of the
nonuniform $(h,k,\mu,\nu)$-dichotomy.
\begin{example}\label{exaaa}
Consider the differential equation in $\R^2$
\begin{equation}\label{eqxy}
z_1'=(-\eta_1\hat{h}'(t)/\hat{h}(t)+\zeta_1(t))z_1,~~~~
z_2'=(\eta_3\hat{k}'(t)/\hat{k}(t)+\zeta_2(t))z_2
\end{equation}
for $t\in \R$, where
$$
\ba{c}
\zeta_1(t)=\eta_2(\hat{\mu}'(t)/\hat{\mu}(t))(\log\hat{\mu}(t)\cos(\log\hat{\mu(t)})-1),\\
\zeta_2(t)=\eta_2(\hat{\nu}'(t)/\hat{\nu}(t))(\log\hat{\nu}(t)\cos(\log\hat{\nu}(t))-1),
\ea
$$
$\hat{h},\hat{k},\hat{\mu},\hat{\nu}$ are growth rates and
$\eta_1,\eta_2,\eta_3$ are positive constants.
\end{example}
Set $P(t)(z_1,z_2)=z_1$ and $Q(t)(z_1,z_2)=z_2$ for $t\in\R$. Then
we have
$$
T(t,s)P(s)=\left(\hat{h}(t)/\hat{h}(s)\right)^{-\eta_1}e^{\eta_2d_1(t)},~~
T(t,s)Q(s)=\left(\hat{k}(t)/\hat{k}(s)\right)^{\eta_3}e^{\eta_2d_2(t)},
$$
where
$$
\ba{l}
d_1(t)=\log\hat{\mu}(t)(\sin\log\hat{\mu}(t)-1)
+\cos\log\hat{\mu}(t)-\cos\log\hat{\mu}(s)-\log\hat{\mu}(s)(\sin\log\hat{\mu}(s)-1),\\
d_2(t)=\log\hat{\nu}(t)(\sin\log\hat{\nu}(t)-1)
+\cos\log\hat{\nu}(t)-\cos\log\hat{\nu}(s)-\log\hat{\nu}(s)(\sin\log\hat{\nu}(s)-1).
\ea
$$
It follows that
$$
\ba{l}
\|T(t,s)P(s)\| \le e^{2 \eta_2}
\left(\hat{h}(t)/\hat{h}(s)\right)^{-\eta_1}
\hat{\mu}(s)^{2 \eta_2}
\leq e^{2 \eta_2}
\left(\hat{h}(t)/\hat{h}(s)\right)^{-\eta_1}
\hat{\mu}(|s|)^{2 \eta_2},~t\geq s,\\
\|T(t,s)Q(s)\|\le e^{2\eta_2}\left(\hat{k}(s)/\hat{k}(t)\right)^{-\eta_3}\hat{\nu}(s)^{2\eta_2}\leq e^{2\eta_2}\left(\hat{k}(s)/\hat{k}(t)\right)^{-\eta_3}\hat{\nu}(|s|)^{2\eta_2},~s\geq t.
\ea
$$
This implies that \eqref{eqxy} admits a nonuniform
$(h,k,\mu,\nu)$-dichotomy with
\[
K=e^{2\eta_2}, \quad a=-\eta_1, \quad b=\eta_3\quad\text{and}
\quad \ve=2\eta_2.
\]
Particularly, when $\hat{h},\hat{k},\hat{\mu},\hat{\nu}$ are
chosen as different functions, we obtain new nonuniform
dichotomies different to the existing ones. For example, if $\hat{h}(t)=t+1$, $\hat{k}(t)=e^t$, $\hat{\mu}(t)=t^2+1$ and
$\hat{\nu}(t)=e^{t^2}$ for $t\in\R^+$, then \eqref{eqxy} admits a dichotomy that can not be covered by any known dichotomy in the literatures.

\section{Existence of nonuniform $(h,k,\mu,\nu)$-dichotomy}\label{existencesection}
\noindent

In this section, in terms of appropriate
Lyapunov exponents and Lyapunov functions, some sufficient criteria are established for linear
dynamical systems in a finite-dimensional space to have a
nonuniform $(h,k,\mu,\nu)$-dichotomy. Those results show that the notion of nonuniform
$(h,k,\mu,\nu)$-dichotomy occurs in a very natural way for linear
nonautonomous dynamical systems.

\subsection{Lyapunov exponents and nonuniform $(h,k,\mu,\nu)$-dichotomy}

Assume that $A(t)$ in \eqref{eqlinear} is a
continuous  $n\times n$ matrix function of block form, i.e., $A(t)={\rm diag}(W_1(t), W_2(t))$ for $t\in\R^+$ and $\R^n=E\oplus F$ , where $\dim E=l$ and $\dim
F=n-l$. For $t\geq 0$, consider
\begin{equation}\label{existenceeqa}
x'_1=W_1(t)x_1,
\end{equation}
\begin{equation}\label{existenceeqb}
x'_2=W_2(t)x_2
\end{equation}
and the corresponding adjoint systems
\begin{equation}\label{existenceeqc}
y'_1=-W_1(t)^*y_1,
\end{equation}
\begin{equation}\label{existenceeqd}
y'_2=-W_2(t)^* y_2
\end{equation}
where $W_1(t)^*$ and $W_2(t)^*$ are the transpose of $W_1(t)$ and $W_2(t)$, respectively. Define $\vp:E\ra[-\iy,+\iy]$
and $\psi:F\ra[-\iy,+\iy]$ by
\begin{equation}\label{eqvpa}
\vp(x_1^0)=\limsup\limits_{t\ra+\iy}\f{\log\|x_1(t)\|}{\log
h(t)}~~\mbox{and}~~\psi(x_2^0)=\limsup\limits_{t\ra+\iy}\f{\log\|x_2(t)\|}{\log
k(t)},
\end{equation}
where $x_1(t)$ is the solution of \eqref{existenceeqa} with
$x_1(0)=x_1^0$ and $x_2(t)$ is the solution of
\eqref{existenceeqb} with $x_2(0)=x_2^0$ (we assume that
$\log0=-\iy$).
Then, one has the following claims
\begin{itemize}
\item [(1)] $\vp(0)=-\iy$ and $\psi(0)=-\iy$;

\item [(2)] $\vp(c x_1^0)=\vp(x_1^0)$ and $\psi(c
x_2^0)=\psi(x_2^0)$ for each  $x_1^0\in E, x_2^0\in F$ and
$c\in\R\setminus\{0\}$;

\item [(3)] for any $x'_1,x''_1\in E$ and $x'_2,x''_2\in F$,
$$\vp(x'_1+x''_1)\leq\max\{\vp(x'_1),\vp(x''_1)\},~~\psi(x'_2+x''_2)\leq\max\{\psi(x'_2),\psi(x''_2)\};$$

\item [(4)] $x_1^1,\cdots,x_1^{m}$ are linearly independent if
$\vp(x_1^1),\cdots,\vp(x_1^m)$ are distinct for
$x_1^1,\cdots,x_1^{m}\in E\setminus\{0\}$;
$x_2^1,\cdots,x_2^{m'}$ are linearly independent if
$\psi(x_2^1),\cdots,\psi(x_2^{m'})$ are distinct for
$x_2^1,\cdots,x_2^{m'}\in F\setminus\{0\}$;

\item [(5)] $\vp$ has at most $r\leq l$ distinct values in
$E\setminus\{0\}$, say $-\iy\leq \la_1<\cdots<\la_r\leq+\iy$; $\psi$ has at most $r'\leq n-l$ distinct values in
$F\setminus\{0\}$, say $ -\iy\leq \chi_1<\cdots<\chi_{r'}\leq+\iy.
$
\end{itemize}
Therefore, from  (1)-(3), it follows that $(\vp,\psi)$ is the so-called $(h,k)$
Lyapunov exponent with respect to the linear equation
\eqref{eqlinear}.

Let  $y_1(t)$ be the solution of \eqref{existenceeqc} with
$y_1(0)=y_1^0$ and $y_2(t)$ be the solution of
\eqref{existenceeqd} with $y_2(0)=y_2^0$. Consider
$\bar{\vp}:E\ra[-\iy,+\iy]$ and $\bar{\psi}:F\ra[-\iy,+\iy]$
defined by
\begin{equation}\label{eqvpb}
\bar{\vp}(y_1^0)=\limsup\limits_{t\ra+\iy}\f{\log\|y_1(t)\|}{\log
\bar{h}(t)}~~\mbox{and}~~\bar{\psi}(y_2^0)=\limsup\limits_{t\ra+\iy}\f{\log\|y_2(t)\|}{\log
\bar{k}(t)},
\end{equation}
where $\bar{h}(t),\bar{k}(t)$ are growth rates.  Then
\begin{itemize}
\item [(6)] $(\bar{\vp},\bar{\psi})$ is the $(\bar{h},\bar{k})$ Lyapunov exponent;

\item [(7)] $\bar{\vp}$ takes at most $\bar{r}\leq l$ distinct
values in $E\setminus\{0\}$, say $ -\iy\leq
\bar{\la}_{\bar{r}}<\cdots<\bar{\la}_1\leq+\iy $; $\bar{\psi}$
takes at most $\bar{r}'\leq n-l$ distinct values in
$F\setminus\{0\}$, say $ -\iy\leq
\bar{\chi}_{\bar{r}'}<\cdots<\bar{\chi}_1\leq+\iy. $
\end{itemize}

Let $\varrho_1,\cdots,\varrho_n$ and $\zeta_1,\cdots,\zeta_n$ be
two bases of $\R^n$, they are said to be dual if
$(\varrho_i,\zeta_j)=\og_{ij}$ for every $i,j$, where
$(\cdot,\cdot)$ is the standard inner product in $\R^n$ and
$\og_{ij}$ is the Kronecker symbol. In order to introduce the
regularity coefficients of $\vp,\bar{\vp}$ and $\psi,\bar{\psi}$,
$\la_i,\bar{\la}_i,\chi_i,\bar{\chi}_i$ are assumed to be finite.

\begin{definition}
The regularity coefficients of $\vp$ and $\bar{\vp}$ is defined by
$$
\ga(\vp,\bar{\vp})=\min\max\{\vp(\da_i)+\bar{\vp}(\bar{\da}_i):1\leq
i\leq l\},
$$
where the minimum is taken over all dual bases
$\da_1,\cdots,\da_l$ and $\bar{\da}_1,\cdots,\bar{\da}_l$ of $E$.
\end{definition}

\begin{definition}
The regularity coefficients of $\psi$ and $\bar{\psi}$ is defined
by
$$
\bar{\ga}(\psi,\bar{\psi})=\min\max\{\psi(\ep_i)+\bar{\psi}(\bar{\ep}_i):1\leq
i\leq n-l\},
$$
where the minimum is taken over all dual bases
$\ep_1,\cdots,\ep_{n-l}$ and $\bar{\ep}_1,\cdots,\bar{\ep}_{n-l}$
of $F$.
\end{definition}


\begin{theorem}
Assume that $\vp(x_1)<0$ for any $x_1\in E\setminus\{0\}$ and
$\psi(x_2)>0$ for any $x_2\in F\setminus\{0\}$ with
$\la_r<0<\chi_1$. Then for any sufficiently small $\tilde{\ve}>0$,
\eqref{eqlinear} with A(t) as in \eqref{eqgghhaa} admits a
nonuniform $(h,k,\mu,\nu)$-dichotomy on $\R^+$ with
$$a=\la_r+\tilde{\ve},~~b=\chi_1+\tilde{\ve},~~\ve=\max\{\ga(\vp,\bar{\vp}),~~\bar{\ga}(\psi,\bar{\psi})\}+\tilde{\ve},~~\mu(t)=h(t)\bar{h}(t),~~\nu(t)=k(t)\bar{k}(t).
$$
\end{theorem}
\begin{proof}
Let $X_1(t)$ be a fundamental solution matrix  of
\eqref{existenceeqa}. It is not difficult to show that
$Y_1(t)=(X_1(t)^*)^{-1}$ is a fundamental solution matrix of
\eqref{existenceeqc}. Let $m_j=\vp(x^j_1(0))$ and
$n_j=\bar{\vp}(y^j_1(0))$ for $j=1,\cdots,l$, where
$x^1_1(t),\cdots,x^l_1(t)$ are the columns of $X_1(t)$ and
$y^1_1(t),\cdots,y^l_1(t)$ are the columns of $Y_1(t)$. For any
$\tilde{\ve}>0$ and $\bar{t}>0$, it follows from  \eqref{eqvpa}
and \eqref{eqvpb} that there exists a constant $\bar{K}_1^1$ such
that $\|x^j_1(t)\|\leq \bar{K}_1^1h(t)^{m_j+\tilde{\ve}}$ and
$\|y^j_1(t)\|\leq \bar{K}_1^1\bar{h}(t)^{n_j+\tilde{\ve}}$ for
$t\geq\bar{t}$ and $j=1,\cdots,l$. On the other hand, there exists
a sufficiently large $\bar{K}_1^2$ such that $\|x^j_1(t)\|\leq
\bar{K}_1^2h(t)^{m_j+\tilde{\ve}}$ and $\|y^j_1(t)\|\leq
\bar{K}_1^2\bar{h}(t)^{n_j+\tilde{\ve}}$ for any $t\in[0,\bar{t}]$
and $j=1,\cdots,l$. Let
$\bar{K}_1=\max\{\bar{K}_1^1,\bar{K}_1^2\}$, then
\begin{equation}\label{existenceproofa}
\|x^j_1(t)\|\leq
\bar{K}_1h(t)^{m_j+\tilde{\ve}},~~\|y^j_1(t)\|\leq
\bar{K}_1\bar{h}(t)^{n_j+\tilde{\ve}},~~t\geq 0,~~j=1,\cdots,l.
\end{equation}
Note that
$Y_1(t)^*X_1(t)=\id$, then $(x_1^i(t),y_1^j(t))=\og_{ij}$ for $i,j=1,\cdots,l$. It is clear that, if the matrix $X_1(t)$ is appropriately select, then
$$\ga(\vp,\bar{\vp})=\max\{m_j+n_j:j=1,\cdots,l\}.$$
Let $U(t,s):=X_1(t)X_1^{-1}(s)$ for $t\geq s$, then
$U(t,s)=X_1(t)Y_1(s)^*$ and the entries of $U(t,s)$ are
$u_{ik}(t,s)=\sum\limits_{j=1}^l x_1^{ij}(t)y_1^{kj}(s)$. It
follows from \eqref{existenceproofa} that
\begin{align*}
|u_{ik}(t,s)|&\leq\sum\limits_{j=1}^l
|x_1^{ij}(t)||y_1^{kj}(s)|\leq\sum\limits_{j=1}^l
\|x_1^{j}(t)\|\|y_1^{j}(s)\|\\
&\leq\sum\limits_{j=1}^l\bar{K}_1^2h(t)^{m_j+\tilde{\ve}}\bar{h}(s)^{n_j+\tilde{\ve}}\\
&\leq
\sum\limits_{j=1}^l\bar{K}_1^2(h(t)/h(s))^{m_j+\tilde{\ve}}h(s)^{m_j+\tilde{\ve}}\bar{h}(s)^{n_j+\tilde{\ve}}\\
&\leq\bar{K}_1^2
l(h(t)/h(s))^{\la_r+\tilde{\ve}}(h(s)\bar{h}(s))^{\ga(\vp,\bar{\vp})+\tilde{\ve}}.
\end{align*}
Let $\xi=\sum\limits_{k=1}^ll_ke_k$ with
$\|\xi\|^2=\sum\limits_{k=1}^ll_k^2=1$, where $e_1,\cdots,e_l$ are
the standard orthogonal basis of $E$. Therefore,
\begin{align*}
\|U(t,s)\xi\|^2&=\left\|\sum\limits_{i=1}^l\sum\limits_{k=1}^ll_ku_{ik}(t,s)e_i\right\|^2
\leq\sum\limits_{i=1}^l\left(\sum\limits_{k=1}^ll_k^2\sum\limits_{k=1}^lu_{ik}(t,s)^2\right)
\leq\sum\limits_{i=1}^l\sum\limits_{k=1}^lu_{ik}(t,s)^2,
\end{align*}
which implies that
\begin{align*}
\|U(t,s)\|&\leq\left(\sum\limits_{i=1}^l\sum\limits_{k=1}^lu_{ik}(t,s)^2\right)^{1/2}\\
&\leq\bar{K}_1^2l^2(h(t)/h(s))^{\la_r+\tilde{\ve}}(h(s)\bar{h}(s))^{\ga(\vp,\bar{\vp})+\tilde{\ve}}\\
&\leq\bar{K}_1^2l^2(h(t)/h(s))^{a}\mu(s)^{\ve}.
\end{align*}

Let $X_2(t)$ be a fundamental solution matrix  of
\eqref{existenceeqb}, then $Y_2(t)=(X_2(t)^*)^{-1}$ is a
fundamental solution matrix of \eqref{existenceeqd}. Let
$\bar{m}_j=\psi(x^j_2(0))$ and $\bar{n}_j=\bar{\psi}(y^j_2(0))$
for $j=1,\cdots,n-l$, where $x^1_2(t),\cdots,x^{n-l}_2(t)$ are the
columns of $X_2(t)$ and $y^1_2(t),\cdots,y^{n-l}_2(t)$ are the
columns of $Y_2(t)$. Proceeding similarly to the above, there
exists a positive constant $\bar{K}_2$ such that
\begin{equation}\label{existenceproofb}
\|x^j_2(t)\|\leq
\bar{K}_2k(t)^{\bar{m}_j+\tilde{\ve}}~~\mbox{and}~~\|y^j_2(t)\|\leq
\bar{K}_2\bar{k}(t)^{\bar{n}_j+\tilde{\ve}},~~t\geq 0,~~j=1,\cdots,n-l.
\end{equation}
Moreover, $(x_2^i(t),y_2^j(t))=\og_{ij}$ for $i,j=1,\cdots,n-l$ since
$Y_2(t)^*X_2(t)=\id$ and $X_2(t)$ can be appropriately
selected such that
$$\bar{\ga}(\psi,\bar{\psi})=\max\{\bar{m}_j+\bar{n}_j:j=1,\cdots,n-l\}.$$
Let $V(t,s)=X_2(t)X_2^{-1}(s)$ for $0\leq t\leq s$.  Then
$V(t,s)=X_2(t)Y_2(s)^*$ and the entries of $V(t,s)$ are
$v_{ik}(t,s)=\sum\limits_{j=1}^{n-l} x_2^{ij}(t)y_2^{kj}(s)$. By
\eqref{existenceproofb}, one has
\begin{align*}
|v_{ik}(t,s)|&\leq\sum\limits_{j=1}^{n-l}
|x_2^{ij}(t)||y_2^{kj}(s)|\leq\sum\limits_{j=1}^{n-l}
\|x_2^{j}(t)\|\|y_2^{j}(s)\|\\
&\leq\sum\limits_{j=1}^{n-l}\bar{K}_2^2k(t)^{\bar{m}_j+\tilde{\ve}}\bar{k}(s)^{\bar{n}_j+\tilde{\ve}}\\
&\leq
\sum\limits_{j=1}^{n-l}\bar{K}_2^2(k(s)/k(t))^{-(\bar{m}_j+\tilde{\ve})}k(s)^{\bar{m}_j+\tilde{\ve}}\bar{k}(s)^{\bar{n}_j+\tilde{\ve}}\\
&\leq\bar{K}_2^2
(n-l)(k(s)/k(t))^{-(\chi_1+\tilde{\ve})}(k(s)\bar{k}(s))^{\bar{\ga}(\psi,\bar{\psi})+\tilde{\ve}}
\end{align*}
and
\begin{align*}
\|V(t,s)\bar{\xi}\|^2&=\left\|\sum\limits_{i=1}^{n-l}\sum\limits_{k=1}^{n-l}\bar{l}_kv_{ik}(t,s)\bar{e}_i\right\|^2
\leq\sum\limits_{i=1}^{n-l}\left(\sum\limits_{k=1}^{n-l}\bar{l}_k^2\sum\limits_{k=1}^{n-l}v_{ik}(t,s)^2\right)
\leq\sum\limits_{i=1}^{n-l}\sum\limits_{k=1}^{n-l}v_{ik}(t,s)^2,
\end{align*}
where $\bar{\xi}=\sum\limits_{k=1}^{n-l}\bar{l}_k\bar{e}_k$ and
$\sum\limits_{k=1}^{n-l}\bar{l}_k^2=1$ with $\bar{e}_1,\cdots,\bar{e}_{n-l}$ being the standard orthogonal basis
of $F$. Therefore,
\begin{align*}
\|V(t,s)\|&\leq\left(\sum\limits_{i=1}^{n-l}\sum\limits_{k=1}^{n-l}v_{ik}(t,s)^2\right)^{1/2}\\
&\leq
\bar{K}_2^2(n-l)^2(k(s)/k(t))^{-(\chi_1+\tilde{\ve})}(k(s)\bar{k}(s))^{\bar{\ga}(\psi,\bar{\psi})+\tilde{\ve}}\\
&\leq\bar{K}_2^2(n-l)^2(k(s)/k(t))^{-b}\nu(s)^{\ve}.
\end{align*}
The proof is complete.
\end{proof}

Intuitively, it seems very restrictive that $A(t)$ is assumed to be of block form and $\la_i,\bar{\la}_i,\chi_i,\bar{\chi}_i$ are assumed to be finite. In fact, from the view point of Lyapunov's theory of regularity and ergodic theory, those assumptions are natural,
typical and quite reasonable. For example,
by the Oseledets-Pesin reduction theorem (see Theorem 3.5.5 in
\cite{bareirapesinbook2007}), there exists a coordinate change
(maintaining the values of the Lyapunov exponents) transforming
the matrix $A(t)$ into a block form for $\mu$-almost every $x$
with respect to a time-independent decomposition. We refer the
reader to \cite{bareirapesinbook2007} for a more detailed
exposition.

Although the discussion is carried out for the relatively special
case, in fact, one can confirm the existence of nonuniform
$(h,k,\mu.\nu)$-dichotomy for more general dynamical systems. For
example, consider the linear systems
\begin{equation}\label{eqlinear2}
x'=\wa(t) x,~t\in\R^+
\end{equation}
and
\begin{equation}\label{eqlinear3}
y'=\wb(t)
y=\begin{pmatrix}\wb_1(t)&0\\0&\wb_2(t)\end{pmatrix}y,~t\in\R^+,
\end{equation}
where $\wa(t),\wb(t)$ are continuous  $n\times n$ matrix
functions, $\wb_1(t)$ and $\wb_2(t)$ are matrices of lower order
than $\wb(t)$.  \eqref{eqlinear2} is said to be \emph{reducible}
if there exist a continuously differentiable invertible matrix
$S(t)$ and a constant $\widetilde{M}>0$ such that
$$
S'=\wa(t)S-S\wb(t),~~\|S(t)\|\leq\widetilde{M},~~\|S^{-1}(t)\|\leq\widetilde{M},~~t\in\R^+.$$
Direct calculation shows that, if
$y(t)$ is a solution of \eqref{eqlinear3}, then $x(t)=S(t)y(t)$ is
a solution of \eqref{eqlinear2}. Therefore, if \eqref{eqlinear2}
is reducible and \eqref{eqlinear3} admits a nonuniform
$(h,k,\mu.\nu)$-dichotomy, then \eqref{eqlinear2} also admits a
nonuniform $(h,k,\mu.\nu)$-dichotomy.

\subsection{Lyapunov functions and nonuniform $(h,k,\mu,\nu)$-dichotomy}

Consider the function
\begin{equation}\label{Y-3}
H(t,x)=\langle S(t)x,x\rangle,~~t\in\R,~x\in\R^n,
\end{equation}
where $S(t)$ is a given $n\times n$
matrix function and
\begin{equation}\label{Y-4}
\dot{H}(t,x)=\frac{d}{dh}H(t+h,T(t+h,t)x)|_{h=0}.
\end{equation}
For fixed $\tau\in\R$, set
\begin{equation}\label{eqzxaaa}
\begin{split}
E_\tau^s&:=\{0\}\cup\{x\in \R^n:H(t,T(t,\tau)x)>0,t\geq \tau\},\\
E_\tau^u&:=\{0\}\cup\{x\in \R^n:H(t,T(t,\tau)x)<0,t\geq \tau\}.
\end{split}
\end{equation}
Similar to Lemma 5 in
\cite{Barreira2009bb}, it is not difficult to show that  $E_\tau^s$ and
$E_\tau^u$ are both subspaces and $E_\tau^s\oplus
E_\tau^u=\mathbb{R}^n$ for $\tau\in\R$. Let
$U(t,\tau)=T(t,\tau)|_{E_\tau^s}$ and
$V(t,\tau)=T(t,\tau)|_{E_\tau^u}$ for $t\geq\tau$. Then, one has
\begin{equation}\label{eqzb}
U(t,\tau)(E_\tau^s)=E_t^s~~\mbox{and} ~~V(t,\tau)(E_\tau^u)=E_t^u.
\end{equation}

\begin{lemma}[\cite{Barreira2013}]\label{leqqqaaa}
Given continuous functions $w, f:[\tau, t]\ra \R^+$ and
$\hat{\eta}>0$, if $w(x)-w(\tau)\geq \hat{\eta}\int_\tau^x
w(z)f(z) dz$ for $x\in[\tau, t]$, then $w(x)\geq
w(\tau)\exp(\hat{\eta}\int_\tau^xf(z)dz)$ for $x \in[\tau,
t].$
\end{lemma}

\begin{theorem}\label{thlx}
Assume that
\begin{itemize}
\item[\tu{(i)}] there exists a symmetric invertible $n\times n$
matrices $S\in C^1(\R,\R^{n\times n})$  such that
\begin{equation}\label{jcyaaa}
\limsup_{t\rightarrow\pm\infty}\frac{\log\|S(t)\|}{\log(\mu(|t|)^\ve+\nu(|t|)^\ve)}<\infty
\end{equation}
and
\begin{equation}\label{jcybbb}
S'(t)+S(t)A(t)+A(t)^*S(t)\leq-\id,~~t\in\R;
\end{equation}

\item[\tu{(ii)}] for a given $\tau\in\R$, there exist constants
$\hat{\eta}_i>0, i=1,2$ such that
\begin{equation}\label{jcyccc}
x\in E_\tau^s,~
\dot{H}(t,U(t,\tau)x)\leq-\hat{\eta}_1(h'(t)/h(t))|H(t,U(t,\tau)x)|,
\end{equation}
\begin{equation}\label{jcyddd}
x\in E_\tau^u,~
\dot{H}(t,V(t,\tau)x)\leq-\hat{\eta}_2(k'(t)/k(t))|H(t,V(t,\tau)x)|;
\end{equation}

\item[\tu{(iii)}] there exist constants $\hat{d}>0$,
$\hat{k}_i,\hat{l}_i\geq 0, i=1,2$ such that
\begin{equation}\label{jcyeee}
\|U(t,\tau)\|\leq \hat{l}_1\mu(t)^{\hat{k}_1}, ~|t-\tau|\leq \hat{d}
\end{equation}
and
\begin{equation}\label{jcyfff}
\|V(t,\tau)\|\leq \hat{l}_2\nu(t)^{\hat{k}_2}, ~|t-\tau|\leq \hat{d};
\end{equation}

\item[\tu{(iv)}] $h(t)/h(\tau)\geq\mu(t)/\mu(\tau)$  for $t\geq \tau$ when
$\hat{\eta}_1>2\hat{k}_1$.
\end{itemize}
Then \eqref{eqlinear} has a nonuniform $(h,k,\mu,\nu)$-dichotomy
on $\R$.
\end{theorem}

\begin{proof}
First, we show that, for $\tau\in\R$ and $t\geq\tau$,
\begin{equation}\label{equabc}
|H(\tau,x)|\geq
\f{\hat{d}}{\hat{l}^2_1}\mu(\tau)^{-2\hat{k}_1}\|x\|^2,~x\in
E_\tau^s;~~ |H(\tau,x)|\geq
\f{\hat{d}}{\hat{l}^2_2}\nu(\tau)^{-2\hat{k}_2}\|x\|^2,~x\in
E_\tau^u,
\end{equation}
and
\begin{equation}\label{eq12jcy}
\begin{split}
H(t,U(t,\tau)x)&\leq\left(\f{h(t)}{h(\tau)}\right)^{-\hat{\eta}_1}H(\tau,x),~~x\in E_\tau^s,\\
|H(t,V(t,\tau)x)|&\geq\left(\f{k(t)}{k(\tau)}\right)^{\hat{\eta}_2}|H(\tau,x)|,~~x\in
E_\tau^u.
\end{split}
\end{equation}
In fact, let $x(t)=U(t,\tau)x$ for $x\in E_\tau^s$. By \eqref{Y-3} and \eqref{jcybbb}, one has
\begin{equation}\label{Z-18}
\begin{split}
\f{d}{dt}H(t,x(t))&=\langle S'(t)x(t),x(t)\rangle+\langle S(t)x'(t),x(t)\rangle+\langle S(t)x(t),x'(t)\rangle\\
&=
\langle(S'(t)+S(t)A(t)+A(t)^*S(t))x(t),x(t)\rangle\leq-\|x(t)\|^2.
\end{split}
\end{equation}
It follows from $H(\tau+\hat{d},x(\tau+\hat{d}))\geq0, x\in E_\tau^s$, \eqref{jcyeee} and \eqref{Z-18}
that
\begin{align*}
H(\tau,x)&\geq H(\tau,x)-H(\tau+\hat{d},x(\tau+\hat{d}))\\
&=-\int_\tau^{\tau+\hat{d}}\f{d}{dr}H(r,x(r))dr\geq\int_\tau^{\tau+\hat{d}}\|x(r)\|^2dr\\
&=\int_\tau^{\tau+\hat{d}}\|U(r,\tau)x\|^2dr\geq\|x\|^2\int_\tau^{\tau+\hat{d}}\f{dr}{\|U(\tau,r)\|^2}\\
&\geq\|x\|^2\int_\tau^{\tau+\hat{d}}\f{1}{\hat{l}^2_1}\mu(\tau)^{-2\hat{k}_1}dr
=\f{\hat{d}}{\hat{l}^2_1}\mu(\tau)^{-2\hat{k}_1}\|x\|^2.
\end{align*}
That is, \eqref{equabc} is valid.

For $x\in E_\tau^s$, define
$\theta:[\tau,t]\rightarrow \mathbb{R}^+$ by
$\theta(s)=H(t+\tau-s,U(t+\tau-s,\tau)x)$. \eqref{jcyccc} together with Lemma \ref{leqqqaaa} gives
\begin{align*}
\theta(\tau)-\theta(t)&=H(t,U(t,\tau)x)-H(\tau,x)=\int_\tau^t\dot{H}(v,U(v,\tau)x)dv\leq-\hat{\eta}_1\int_\tau^t\f{h'(v)}{h(v)}H(v,U(v,\tau)x)dv\\
&=
-\hat{\eta}_1\int_\tau^t\f{h'(v)}{h(v)}\theta(t+\tau-v)dv=-\hat{\eta}_1\int_\tau^{t}\frac{h'(t+\tau-s)}{h(t+\tau-s)}\theta(s)ds,
\end{align*}
which means
$$
\theta(t)-\theta(\tau)\geq\hat{\eta}_1\int_\tau^{t}\f{h'(t+\tau-s)}{h(t+\tau-s)}\theta(s)ds
$$
and
$$
\theta(t)\geq\theta(\tau)(h(t)/h(\tau))^{\hat{\eta}_1}.
$$
Then the first inequality of \eqref{eq12jcy} holds.

For $x\in E_\tau^u$, note that
$H(\tau-\hat{d},x(\tau-\hat{d}))\leq0$, by \eqref{jcyfff} and
\eqref{Z-18}, one has
\begin{align*}
|H(\tau,x)|&\geq|H(\tau,x)|-|H(\tau-\hat{d},x(\tau-\hat{d}))|
=H(\tau-\hat{d},x(\tau-\hat{d}))-H(\tau,x)\\
&=-\int_{\tau-\hat{d}}^\tau\f{d}{dr}H(r,x(r))dr\geq\int_{\tau-\hat{d}}^\tau\|x(r)\|^2dr=\int_{\tau-\hat{d}}^\tau\|V(r,\tau)x\|^2dr\\
& \geq\|x\|^2\int_{\tau-\hat{d}}^\tau\f{dr}{\|V(\tau,r)\|^2}dr\geq
\f{\hat{d}}{\hat{l}^2_2}\nu(\tau)^{-2\hat{k}_2}\|x\|^2.
\end{align*}
For a given $\tau\in\R$ and any $x\in E_\tau^u$, it follows from
\eqref{jcyddd} that
\begin{align*}
|H(t,V(t,\tau))x|-|H(\tau,x)|&=H(\tau,V(\tau,\tau)x)-H(t,V(t,\tau)x)\\
&=-\int_\tau^t\dot{H}(v,V(v,\tau)x)dv\geq\hat{\eta}_2\int_\tau^t\f{k'(v)}{k(v)}\left|H(v,V(v,\tau)x)\right|dv,~~t\geq\tau.
\end{align*}
By Lemma \ref{leqqqaaa}, the second inequality of
\eqref{eq12jcy} holds.

By \eqref{jcyaaa}, direct calculation  shows that there exist constants $\hat{a},
\hat{b}$ such that
\begin{equation}\label{eqxwq}
\|S(t)\|\leq\hat{a} (\mu(|t|)^\ve+\nu(|t|)^\ve)^{\hat{b}},~~t\in\R.
\end{equation}

Next we establish the norm
bounds of the evolution operators $U(t,\tau)$ and $V(t,\tau)$,
that is, for $\tau\in\R$ and $t\geq\tau$.
\begin{equation}\label{eq1234}
\begin{split}
\|U(t,\tau)\|^2&\leq\f{\hat{a}\hat{l}^2_1}{\hat{d}}\left(\f{h(t)}{h(\tau)}\right)^{-\hat{\eta}_1+2\hat{k}_1}
(\mu(|\tau|)^\ve+\nu(|\tau|)^\ve)^{\hat{b}}\mu(\tau)^{2\hat{k}_1},~x\in E_\tau^s,\\
\|V(t,\tau)^{-1}\|^2&\leq\f{\hat{a}\hat{l}^2_2}{\hat{d}}\left(\f{k(t)}{k(\tau)}\right)^{-\hat{\eta}_2}(\mu(|t|)^\ve+\nu(|t|)^\ve)^{\hat{b}}\nu(t)^{2\hat{k}_2},~x\in
E_\tau^u.
\end{split}
\end{equation}

For any $x\in
E_\tau^s$ and $t\geq\tau$, it follows from \eqref{equabc},
\eqref{eq12jcy},\eqref{eqzb} and condition (iv) that
\begin{align*}
\|U(t,\tau)x\|^2&\leq\f{\hat{l}^2_1}{\hat{d}}\mu(t)^{2\hat{k}_1}|H(t,U(t,\tau)x)|\leq\f{\hat{l}^2_1}
{\hat{d}}\left(\f{h(t)}{h(\tau)}\right)^{-\hat{\eta}_1}\mu(t)^{2\hat{k}_1}H(\tau,x)\\
&\leq\f{\hat{l}^2_1}
{\hat{d}}\left(\f{h(t)}{h(\tau)}\right)^{-\hat{\eta}_1}\mu(t)^{2\hat{k}_1}\|S(\tau)\|\|x\|^2\\
&\leq\f{\hat{a}\hat{l}^2_1}
{\hat{d}}\left(\f{h(t)}{h(\tau)}\right)^{-\hat{\eta}_1}\mu(t)^{2\hat{k}_1}
(\mu(|\tau|)^\ve+\nu(|\tau|)^\ve)^{\hat{b}}\|x\|^2\\
&=\f{\hat{a}\hat{l}^2_1}
{\hat{d}}\left(\f{h(t)}{h(\tau)}\right)^{-\hat{\eta}_1}\left(\f{\mu(t)}{\mu(\tau)}\right)^{2\hat{k}_1}
(\mu(|\tau|)^\ve+\nu(|\tau|)^\ve)^{\hat{b}}\mu(\tau)^{2\hat{k}_1}\|x\|^2\\
&\leq\f{\hat{a}\hat{l}^2_1}
{\hat{d}}\left(\f{h(t)}{h(\tau)}\right)^{-\hat{\eta}_1+2\hat{k}_1}
(\mu(|\tau|)^\ve+\nu(|\tau|)^\ve)^{\hat{b}}\mu(\tau)^{2\hat{k}_1}\|x\|^2,
\end{align*}
which implies that the first inequality of \eqref{eq1234} holds.

For any $x\in E_\tau^u$, by \eqref{eqxwq} and \eqref{eqzb}, one has
\begin{align*}
|H(t,V(t,\tau)x)|\leq
\hat{a}(\mu(|t|)^\varepsilon+\nu(|t|)^\varepsilon)^{\hat{b}}\|V(t,\tau)x\|^2.
\end{align*}
Then
\begin{align*}
\|V(t,\tau)x\|^2&\geq\f{1}{\hat{a}}(\mu(|t|)^\varepsilon+\nu(|t|)^\varepsilon)^{-\hat{b}}\left(\f{k(t)}{k(\tau)}\right)^{\hat{\eta}_2}
\f{\hat{d}}{\hat{l}_2^2}\nu(\tau)^{-2\hat{k}_2}\|x\|^2\\
&\geq\f{\hat{d}}{\hat{a}\hat{l}^2_2}\nu(t)^{-2\hat{k}_2}(\mu(|t|)^\varepsilon+\nu(|t|)^\varepsilon)^{-\hat{b}}
\left(\f{k(t)}{k(\tau)}\right)^{\hat{\eta}_2}\left(\f{\nu(t)}{\nu(\tau)}\right)^{2\hat{k}_2}\|x\|^2\\
&\geq\f{\hat{d}}{\hat{a}\hat{l}^2_2}\nu(t)^{-2\hat{k}_2}(\mu(|t|)^\varepsilon+\nu(|t|)^\varepsilon)^{-\hat{b}}
\left(\f{k(t)}{k(\tau)}\right)^{\hat{\eta}_2}\|x\|^2,
\end{align*}
that is, the second inequality of \eqref{eq1234} holds.

Let $P(t):\R^n\ra E_t^s$, $Q(t):\R^n\ra E_t^u$, $t\in\R$ be the projections. Then
\begin{equation}\label{eq2233}
\begin{split}
\|T(t,\tau)P(\tau)\|\leq\|T(t,\tau)|_{E_\tau^t}\|\|
P(\tau)\|=\|U(t,\tau)|\|\| P(\tau)\|,~t\geq\tau,\\
\|T(t,\tau)Q(\tau)\|\leq\|T(t,\tau)|_{E_\tau^t}\|\|
Q(\tau)\|=\|V(t,\tau)|\|\| Q(\tau)\|,~t\leq\tau.
\end{split}
\end{equation}
For any $z\in\R^n$, we have
$z=P(t)z+Q(t)z$, $P(t)z\in E_t^s$ and $Q(t)z\in E_t^u$. By
\eqref{equabc}, one has
$$
-H(t,P(t)z)+\f{\hat{d}}{\hat{l}_1^2}\mu(t)^{-2\hat{k}_1}\|P(t)z\|^2+H(t,Q(t)z)+
\f{\hat{d}}{\hat{l}_2^2}\nu(t)^{-2\hat{k}_2}\|Q(t)z\|^2\leq0.
$$
Note that
\begin{align*}
&\f{\hat{d}}{\hat{l}_1^2}\mu(t)^{-2\hat{k}_1}\|P(t)z-\f{\hat{l}_1^2}{2\hat{d}}\mu(t)^{2\hat{k}_1}S(t)z\|^2+
\f{\hat{d}}{\hat{l}_2^2}\nu(t)^{-2\hat{k}_2}\|Q(t)z+\f{\hat{l}_2^2}{2\hat{d}}\nu(t)^{2\hat{k}_2}S(t)z\|^2\\
&=\f{\hat{d}}{\hat{l}_1^2}\mu(t)^{-2\hat{k}_1}\|P(t)z\|^2+\f{\hat{d}}{\hat{l}_2^2}\nu(t)^{-2\hat{k}_2}
\|Q(t)z\|^2-H(t,P(t)z)+H(t,Q(t)z)\\
&\quad+\f{\hat{l}_1^2}{4\hat{d}}\mu(t)^{2\hat{k}_1}\|S(t)z\|^2+\f{\hat{l}_2^2}{4\hat{d}}\nu(t)^{2\hat{k}_2}\|S(t)z\|^2\\
&\leq\f{\hat{l}_1^2}{4\hat{d}}\mu(t)^{2\hat{k}_1}\|S(t)z\|^2+\f{\hat{l}_2^2}{4\hat{d}}\nu(t)^{2\hat{k}_2}\|S(t)z\|^2.
\end{align*}
In order to complete the proof, one only needs to prove norm estimation
of $P(t)$ and $Q(t)$.  The discussion is divided into two cases.

Case 1.~
If $(\hat{d}/\hat{l}_1^2)\mu(t)^{-2\hat{k}_1}\leq
\hat{d}/(\hat{l}_2^2)\nu(t)^{-2\hat{k}_2}$, then
$$
\|P(t)z-\f{\hat{l}_1^2}{2\hat{d}}\mu(t)^{2\hat{k}_1}S(t)z\|^2+\|Q(t)z+\f{\hat{l}_2^2}{2\hat{d}}\nu(t)^{2\hat{k}_2}S(t)z\|^2
\leq\f{\hat{l}_1^4}{2\hat{d}^2}\mu(t)^{4\hat{k}_1}\|S(t)z\|^2.
$$
Whence,
\begin{align*}
\|P(t)z\|&=\|P(t)z-\f{\hat{l}_1^2}{2\hat{d}}\mu(t)^{2\hat{k}_1}S(t)z+\f{\hat{l}_1^2}{2\hat{d}}\mu(t)^{2\hat{k}_1}S(t)z\|\\
&\leq\|P(t)z-\f{\hat{l}_1^2}{2\hat{d}}\mu(t)^{2\hat{k}_1}S(t)z\|+\f{\hat{l}_1^2}{2\hat{d}}\mu(t)^{2\hat{k}_1}\|S(t)z\|\\
&\leq\f{\sqrt{2}\hat{l}_1^2}{2\hat{d}}\mu(t)^{2\hat{k}_1}\|S(t)z\|+\f{\hat{l}_1^2}{2\hat{d}}\mu(t)^{2\hat{k}_1}\|S(t)z\|\\
&=\f{(\sqrt{2}+1)\hat{l}_1^2}{2\hat{d}}\mu(t)^{2\hat{k}_1}\|S(t)z\|.
\end{align*}
and
\begin{align*}
\|Q(t)z\|&=\|Q(t)z+\f{\hat{l}_2^2}{2\hat{d}}\nu(t)^{2\hat{k}_2}S(t)z-\f{\hat{l}_2^2}{2\hat{d}}\nu(t)^{2\hat{k}_2}S(t)z\|\\
&\leq\|Q(t)z+\f{\hat{l}_2^2}{2\hat{d}}\nu(t)^{2\hat{k}_2}S(t)z\|+\f{\hat{l}_2^2}{2\hat{d}}\nu(t)^{2\hat{k}_2}\|S(t)z\|\\
&\leq\f{\sqrt{2}\hat{l}_1^2}{2\hat{d}}\mu(t)^{2\hat{k}_1}\|S(t)z\|+\f{\hat{l}_1^2}{2\hat{d}}\mu(t)^{2\hat{k}_1}\|S(t)z\|\\
&=\f{(\sqrt{2}+1)\hat{l}_1^2}{2\hat{d}}\mu(t)^{2\hat{k}_1}\|S(t)z\|.
\end{align*}

Case 2.~If $(\hat{d}/\hat{l}_1^2)\mu(t)^{-2\hat{k}_1}>
\hat{d}/(\hat{l}_2^2)\nu(t)^{-2\hat{k}_2}$, then
$$
\|P(t)z-\f{\hat{l}_1^2}{2\hat{d}}\mu(t)^{2\hat{k}_1}S(t)z\|^2+\|Q(t)z+\f{\hat{l}_2^2}{2\hat{d}}\nu(t)^{2\hat{k}_2}S(t)z\|^2
\leq\f{\hat{l}_2^4}{2\hat{d}^2}\nu(t)^{4\hat{k}_2}\|S(t)z\|^2.
$$
Similarly, one has
$$
\|P(t)z\|\leq\f{(\sqrt{2}+1)\hat{l}_2^2}{2\hat{d}}\nu(t)^{2\hat{k}_2}\|S(t)z\|,~\|Q(t)z\|\leq\f{(\sqrt{2}+1)\hat{l}_2^2}{2\hat{d}}\nu(t)^{2\hat{k}_2}\|S(t)z\|.
$$
Then
\begin{equation}\label{eqqqaabb}
\begin{split}
\|P(t)\|\leq\max\left\{\f{(\sqrt{2}+1)\hat{l}_1^2}{2\hat{d}}\mu(t)^{2\hat{k}_1},\f{(\sqrt{2}+1)\hat{l}_2^2}{2\hat{d}}\nu(t)^{2\hat{k}_2}\right\}\|S(t)\|,\\
\|Q(t)\|\leq\max\left\{\f{(\sqrt{2}+1)\hat{l}_1^2}{2\hat{d}}\mu(t)^{2\hat{k}_1},\f{(\sqrt{2}+1)\hat{l}_2^2}{2\hat{d}}\nu(t)^{2\hat{k}_2}\right\}\|S(t)\|.
\end{split}
\end{equation}
It follows from \eqref{eqxwq}, \eqref{eq1234},
\eqref{eq2233} and \eqref{eqqqaabb} that Theorem \ref{thlx} is valid.
\end{proof}

In order to further characterize the nonuniform $(h,k,\mu,\nu)$-dichotomy,
we establish a necessary condition for the existence of  nonuniform $(h,k,\mu,\nu)$-dichotomies for the linear system \eqref{eqlinear}.

\begin{theorem}
If \eqref{eqlinear} admits a nonuniform $(h,k,\mu,\nu)$-dichotomy
\eqref{deeqaaa} on $\R$  with $h,k\in C^1(\R,\R^+)$, then there
exists a symmetric invertible matrix function $S\in
C^1(\R,\R^{n\times n})$ such that
\begin{equation}\label{Z-1}
\limsup_{t\rightarrow\pm\infty}\f{\log\|S(t)\|}{\log(\mu(|t|)^\varepsilon+\nu(|t|)^\varepsilon)}<\infty
\end{equation}
and
\begin{equation}\label{Z-2}
S'(t)+S(t)A(t)+A(t)^*S(t)\leq-\left(P(t)^*P(t)\f{h'(t)}{h(t)}+Q(t)^*Q(t)\frac{k'(t)}{k(t)}\right),~~t\in\R.
\end{equation}
Moreover, there exist constants $\bar{k}_1$ and $\bar{k}_2$ such that
\begin{equation}\label{Z-3}
\begin{split}
\dot{H}(t,T(t,\tau)x)&\leq-\bar{k}_1(h'(t)/h(t))|H(t,T(t,\tau)x)|,~
x\in F_\tau^s,\\
\dot{H}(t,T(t,\tau)
x)&\leq-\bar{k}_2(k'(t)/k(t))|H(t,T(t,\tau)x)|,~x\in F_\tau^u
\end{split}
\end{equation}
for given $\tau\in\R$, where $F_\tau^s=P(\tau)(\R^n)$,
$F_\tau^u=Q(\tau)(\R^n)$ and  $H(t,x)$ is as in \eqref{Y-3}.
\end{theorem}
\begin{proof}
For some positive constant $0<\bar{d}<\min\{-a,b\}$, set
\begin{equation}\label{eqzbj111}
\begin{split}
S(t)=&\int_t^\infty T(v,t)^*P(v)^*P(v)T(v,t)\left(\frac{h(v)}{h(t)}\right)^{-2(a+\bar{d}))}\f{h'(v)}{h(v)}dv\\
     &-\int_{-\infty}^tT(v,t)^*Q(v)^*Q(v)T(v,t)\left(\frac{k(t)}{k(v)}\right)^{2(b-\bar{d})}\f{k'(v)}{k(v)}dv.
\end{split}
\end{equation}
It is clear that $S(t)$ is symmetric for each
$t\in\R$. Note that $\partial T(\tau,t)/\partial
t=-T(\tau,t)A(t)$ and $\partial T(\tau,t)^*/\partial
t=-A(t)^*T(\tau,t)^*$, then $S(t)$ is continuously differentiable.
From \eqref{Y-3}, it follows that
\begin{equation}\label{eqxwq123}
\begin{split}
H(t,x)&=\int_t^\infty\|T(v,t)P(t)x\|^2\left(\f{h(v)}{h(t)}\right)^{-2(a+\bar{d})}\f{h'(v)}{h(v)}dv\\
&\quad-\int_{-\infty}^t\|T(v,t)Q(t)x\|^2\left(\f{k(t)}{k(v)}\right)^{2(b-\bar{d})}\f{k'(v)}{k(v)}dv,~~x\in\R^n.
\end{split}
\end{equation}
If \eqref{eqlinear} admits a nonuniform
$(h,k,\mu,\nu)$-dichotomy, then $\mathbb{R}^n=F_t^s\oplus F_t^u$,
where $F_t^s=P(t)(\R^n)$ and $F_t^u=Q(t)(\R^n)$. Note that
$\langle S(t)x,x\rangle=H(t,x)>0$ for $x\in F_t^s\setminus \{0\}$ and
$\langle S(t)x,x\rangle=H(t,x)<0$ for $x\in F_t^u\setminus \{0\}$, then
$S(t)|_{F_t^s}$ and $S(t)|_{F_t^u}$ are both invertible. This
implies that $S(t)$ is invertible for each $t\in\R$.

By \eqref{deeqaaa}, one has
\begin{align*}
L_1&=:\int_t^\infty\|T(v,t)P(t)\|^2\left(\f{h(v)}{h(t)}\right)^{-2(a+\bar{d})}\f{h'(v)}{h(v)}dv\\
&\leq
K^2\int_t^\infty\left(\f{h(v)}{h(t)}\right)^{2a}\mu(|t|)^{2\varepsilon}
\left(\f{h(v)}{h(t)}\right)^{-2(a+\bar{d})}\f{h'(v)}{h(v)}dv\leq (K^2/2\bar{d})\mu(|t|)^{2\ve}
\end{align*}
and
\begin{align*}
L_2&=:\int_{-\infty}^t\|T(v,t)Q(t)\|^2\left(\f{k(t)}{k(v)}\right)^{2(b-\bar{d})}\f{k'(v)}{k(v)}dv\\
&\leq
K^2\int_{-\infty}^t\left(\f{k(t)}{k(v)}\right)^{-2b}v(|t|)^{2\varepsilon}
\left(\f{k(t)}{k(v)}\right)^{2(b-\bar{d})}\f{k'(v)}{k(v)}dv\leq (K^2/2\bar{d})\nu(|t|)^{2\ve}.
\end{align*}
Then
\begin{equation}\label{eqswk111}
\begin{split}
\|S(t)\|&=\sup\limits_{x\neq0}\f{|H(t,x)|}{\|x\|^2}\leq\sup\limits_{x\neq0}\f{\|S(t)x\|\|x\|}{\|x\|^2}\leq L_1+L_2\leq(K^2/2\bar{d})(\mu(|t|)^{2\ve}+\nu(|t|)^{2\ve}),
\end{split}
\end{equation}
which implies that
\begin{align*}
\limsup_{t\rightarrow{\pm\infty}}\f{\log\|S(t)\|}{\log(\mu(|t|)^\varepsilon+\nu(|t|)^\varepsilon)}
&\leq\limsup_{t\rightarrow{\pm\infty}}\f{\log({K^2}/{2\bar{d}})(\mu(|t|)^{2\varepsilon}+\nu(|t|)^{2\varepsilon})}
    {\log(\mu(|t|)^\varepsilon+\nu(|t|)^\varepsilon)}\\
&\leq\limsup_{t\rightarrow{\pm\infty}}\f{\log({K^2}/{2\bar{d}})}
    {\log(\mu(|t|)^\varepsilon+\nu(|t|)^\varepsilon)}+2<{+\infty}.
\end{align*}
Direct calculation leads to
\begin{equation}\label{eq123}
\begin{split}
S'(t)=&-P(t)^*P(t)\f{h'(t)}{h(t)}-Q(t)^*Q(t)\f{k'(t)}{k(t)}-A(t)^*S(t)-S(t)A(t)\\
      &+2(a+\bar{d})\f{h'(t)}{h(t)}\int_t^\infty T(v,t)^*P(v)^*P(v)T(v,t)\left(\f{h(v)}{h(t)}\right)^{-2(a+\bar{d})}
       \f{h'(v)}{h(v)}dv\\
      &-2(b-\bar{d})\f{k'(t)}{k(t)}\int_{-\infty}^t T(v,t)^*Q(v)^*Q(v)T(v,t)\left(\f{k(t)}{k(v)}\right)^{2(b-\bar{d})}
       \f{k'(v)}{k(v)}dv.
\end{split}
\end{equation}
Thus, \eqref{Z-2} holds since $\bar{d}<\min\{-a,b\}$.

For each given $\tau\in\R$ and any $x\in \R^n$, by \eqref{eq123} and $h'(t),k'(t)\geq 0$, we have
\begin{align*}
\f{dH(t,T(t,\tau)x)}{dt}&=\langle S'(t)T(t,\tau)x,T(t,\tau)x\rangle+\langle S(t)\f{\partial T(t,\tau)}{\partial t}x,T(t,\tau)x\rangle
+\langle S(t)T(t,\tau)x,\f{\partial T(t,\tau)}{\partial  t}x\rangle\\
&=\langle(S'(t)+S(t)A(t)+A(t)^*S(t))T(t,\tau)x,T(t,\tau)x\rangle\\
&\leq -\left\langle\left(P(t)^*P(t)\f{h'(t)}{h(t)}+Q(t)^*Q(t)\f{k'(t)}{k(t)}\right)T(t,\tau)x,T(t,\tau)x\right\rangle\\
&\quad+2(a+\bar{d})\f{h'(t)}{h(t)}\int_t^\infty \|T(v,t)P(t)T(t,\tau)x\|^2\left(\f{h(v)}{h(t)}\right)^{-2(a+\bar{d})}\f{h'(v)}{h(v)}dv\\
&\quad-2(b-\bar{d})\f{k'(t)}{k(t)}\int_{-\infty}^t
 \|T(v,t)Q(t)T(t,\tau)x\|^2\left(\f{k(t)}{k(v)}\right)^{2(b-\bar{d})}\f{k'(v)}{k(v)}dv\\
&\leq2(a+\bar{d})\f{h'(t)}{h(t)}\int_t^\infty \|T(v,\tau)P(\tau)x\|^2\left(\f{h(v)}{h(t)}\right)^{-2(a+\bar{d})}\f{h'(v)}{h(v)}dv\\
&\quad-2(b-\bar{d})\f{k'(t)}{k(t)}\int_{-\infty}^t
 \|T(v,\tau)Q(\tau)x\|^2\left(\f{k(t)}{k(v)}\right)^{2(b-\bar{d})}\f{k'(v)}{k(v)}dv.
\end{align*}
Moreover, by \eqref{eqxwq123}, one has
\begin{align*}
\f{dH(t,T(t,\tau)x)}{dt}&\leq2(a+\bar{d})\f{h'(t)}{h(t)}\int_t^\infty \|T(v,\tau)P(\tau)x\|^2\left(\f{h(v)}{h(t)}\right)^{-2(a+\bar{d})}\f{h'(v)}{h(v)}dv\\
&=2(a+\bar{d})\f{h'(t)}{h(t)}|H(t,T(t,\tau)x)|,~~x\in F_\tau^s
\end{align*}
\begin{align*}
\f{dH(t,T(t,\tau)x)}{dt}&\leq-2(b-\bar{d})\f{k'(t)}{k(t)}\int_{-\infty}^t
 \|T(v,\tau)Q(\tau)x\|^2\left(\f{k(t)}{k(v)}\right)^{2(b-\bar{d})}\f{k'(v)}{k(v)}dv\\
&=-2(b-\bar{d})\f{k'(t)}{k(t)}|H(t,T(t,\tau)x)|,~~x\in F_\tau^u.
\end{align*}
Therefore, \eqref{Z-3} holds with $\bar{k}_1=-2(a+\bar{d})$ and $\bar{k}_2=2(b-\bar{d})$.
\end{proof}
\begin{corollary}
If (iii) and (iv) in Theorem \ref{thlx} hold and
$$P(t)^*P(t)\f{h'(t)}{h(t)}+Q(t)^*Q(t)\f{k'(t)}{k(t)}\geq \id,$$
then \eqref{eqlinear} admits a nonuniform
$(h,k,\mu,\nu)$-dichotomy as in \eqref{deeqaaa} on $\R$ if and
only if (i) and (ii) in Theorem \ref{thlx} hold.
\end{corollary}

\section{Robustness of nonuniform $(h,k,\mu,\nu)$-dichotomy}\label{robustnesssection}
\noindent

This section focuses on the robustness or roughness of
nonuniform $(h,k,\mu,\nu)$-dichotomy, which is one of the most
important properties of dichotomy. The principal aim is to show
that the nonuniform $(h,k,\mu,\nu)$-dichotomy persists under
sufficiently small linear perturbations of the original dynamics.

We first establish the robustness of nonuniform $(h,k,\mu,\nu)$-dichotomy in a finite-
dimensional space by using  Lyapunov functions. Consider the linear perturbed system
\begin{equation}\label{eqpllinearaaa}
x'=(A(t)+B(t))x,
\end{equation}
where $A,B\in C(\R,\R^{n\times n})$. Let $\htt(t,\tau)$ be the evolution operator associated to \eqref{eqpllinearaaa} and
\begin{equation}\label{eqzxaaa123}
\begin{split}
\widehat{E}_\tau^s&:=\{0\}\cup\{x\in \R^n:H(t,\htt(t,\tau)x)>0,t\geq \tau\},\\
\widehat{E}_\tau^u&:=\{0\}\cup\{x\in
\R^n:H(t,\htt(t,\tau)x)<0,t\geq \tau\}
\end{split}
\end{equation}
for each given $\tau\in\R$. It is not difficult to show that
$\widehat{E}_\tau^s$ and $\widehat{E}_\tau^u$ are both subspaces
and $\widehat{E}_\tau^s\oplus \widehat{E}_\tau^u=\mathbb{R}^n$ for $\tau\in\R$.

\begin{theorem}\label{thxwq}
Assume that \eqref{eqlinear} admits a nonuniform
$(h,k,\mu,\nu)$-dichotomy as in \eqref{deeqaaa} on $\R$ and there
exist positive constants $\hat{l}$, $\hat{\da}$ and $\hat{d}$ such
that
\begin{equation}\label{eqzx111}
\|T(t,\tau)\|\leq\hat{l}\min\{\mu(t)^{2\ve},\nu(t)^{2\ve}\},~~|t-\tau|\leq\hat{d},
\end{equation}
\begin{equation}\label{eqzx222}
\|B(t)\|\leq \hat{\da}(\mu(|t|)+\nu(|t|))^{-2\ve}.
\end{equation}
Moreover, if (iv) in Theorem \ref{thlx} holds and
\begin{equation}\label{eq33445566}
P^*(t)P(t)(h'(t)/h(t))+Q^*(t)Q(t)(k'(t)/k(t))-(\hat{\da}K^2/\bar{d})\id\geq\id,
\end{equation}
for a positive constant $\bar{d}<\min\{-a,b\}$, then
\eqref{eqpllinearaaa} also admits a nonuniform
$(h,k,\mu,\nu)$-dichotomy on $\R$.
\end{theorem}
\begin{proof}
We first claim that condition (iii) in Theorem \ref{thlx} holds.
In fact, by the variation of constant method, one has
$$
\htt(t,\tau)=T(t,\tau)+\il_\tau^t T(t,r)B(r)\htt(r,\tau)dr
$$
for each $t,\tau\in\R$ and $|t-\tau|\leq\bar{d}$. Then
\begin{align*}
\widehat{U}(t,\tau)&=\htt(t,\tau)|_{\widehat{E}_\tau^s}=T(t,\tau)|_{\widehat{E}_\tau^s}+\il_\tau^t T(t,r)B(r)\htt(r,\tau)|_{\widehat{E}_\tau^s}dr\\
&=T(t,\tau)|_{\widehat{E}_\tau^s}+\il_\tau^t
T(t,r)B(r)\widehat{U}(r,\tau)dr.
\end{align*}
By \eqref{eqzx111} and \eqref{eqzx222}, we have
\begin{align*}
\|\widehat{U}(t,\tau)\|&\leq\hat{l}\mu(t)^{2\ve}+\hat{l}\hat{\da}\mu(t)^{2\ve}\il_\tau^t(\mu(|r|)+\nu(|r|))^{-2\ve}\|\widehat{U}(r,\tau)\|dr\\
&\leq\hat{l}\mu(t)^{2\ve}+\hat{l}\hat{\da}\mu(t)^{2\ve}\il_\tau^t\mu(r)^{-2\ve}\|\widehat{U}(r,\tau)\|dr.
\end{align*}
From the Gronwall's inequality, it follows that, for $|t-\tau|\leq\hat{d}$,
$$
\|\widehat{U}(t,\tau)\|\leq\hat{l}\exp\{\hat{l}\hat{\da}\hat{d}\}\mu(t)^{2\ve},~~
\|\widehat{V}(t,\tau)\|=\|\htt(t,\tau)|_{\widehat{E}_\tau^u}\|\leq\hat{l}\exp\{\hat{l}\hat{\da}\hat{d}\}\nu(t)^{2\ve}.
$$

In order to prove (i) of Theorem \ref{thlx},  consider
the matrix $S(t)$ in \eqref{eqzbj111}. By \eqref{eqswk111},
\eqref{Z-2} and \eqref{eq33445566}, we have
\begin{align*}
&S'(t)+S(t)(A(t)+B(t))+(A^*(t)+B^*(t))S(t)\\
&\leq-P(t)^*P(t)\f{h'(t)}{h(t)}-Q(t)^*Q(t)\f{k'(t)}{k(t)}+2\|S(t)\|\|B^*(t)\|\id\\
&\leq-P(t)^*P(t)\f{h'(t)}{h(t)}-Q(t)^*Q(t)\f{k'(t)}{k(t)}+(\hat{\da}K^2/\bar{d})\id\leq-\id,~~t\in\R,
\end{align*}
which, together with \eqref{Z-1}, implies (i) of Theorem \ref{thlx}.

Direct calculation leads to
\begin{align*}
\f{dH(t,\htt(t,\tau)x)}{dt}=&\langle
(S'(t)+S(t)(A(t)+B(t))+(A^*(t)+B^*(t))S(t))\htt(t,\tau)x,\htt(t,\tau)x\rangle\\
&\leq2(a+\bar{d})\f{h'(t)}{h(t)}\int_t^\infty \|T(v,t)P(t)\htt(t,\tau)x\|^2\left(\f{h(v)}{h(t)}\right)^{-2(a+\bar{d})}\f{h'(v)}{h(v)}dv\\
&\quad-2(b-\bar{d})\f{k'(t)}{k(t)}\int_{-\infty}^t
 \|T(v,t)Q(t)\htt(t,\tau)x\|^2\left(\f{k(t)}{k(v)}\right)^{2(b-\bar{d})}\f{k'(v)}{k(v)}dv
\end{align*}
for each given $\tau\in\R$ and any $x\in \R^n$ since \eqref{eq123}
holds. By \eqref{eqzxaaa123} and \eqref{eqxwq123}, we get
$$
\f{dH(t,\htt(t,\tau)x)}{dt}\leq2(a+\bar{d})\f{h'(t)}{h(t)}|H(t,\htt(t,\tau)x)|,~~x\in\widehat{E}_\tau^s,
$$
$$
\f{dH(t,\htt(t,\tau)x)}{dt}\leq-2(b-\bar{d})\f{h'(t)}{h(t)}|H(t,\htt(t,\tau)x)|,~~x\in\widehat{E}_\tau^u.
$$
The proof is complete.
\end{proof}

Next we characterize robustness of the nonuniform
$(h,k,\mu,\nu)$-dichotomy in the infinite-dimensional space.

Let $X$ be a Banach space and $Y= (Y,|\cdot|)$ be an open subset of the parameter space (also Banach space).
Consider the linear perturbed system with
parameters
\begin{equation}\label{eqpllinear}
x'=(A(t)+B(t,\la))x,
\end{equation}
where $B:\R\times Y\ra\B(X)$.

\begin{theorem}\label{thaaa}
Assume that
\begin{itemize}

\item[\tu{(a$_1$)}]
\eqref{eqlinear} admits a nonuniform $(h,k,\mu,\nu)$-dichotomy as
in \eqref{deeqaaa} on $\R$ with
$\lim_{t\ra\iy}k(t)^{-b}\nu(|t|)^\ve=0$ and
$\lim_{t\ra-\iy}h(t)^{-a}\mu(|t|)^\ve=0$;

\item[\tu{(a$_2$)}] there is a positive constant $N$
such that
$$
\nu(|t|)^\ve\il_{-\iy}^t\mu(|\tau|)^{-\og}d\tau+
\mu(|t|)^\ve\il_t^\iy\nu(|\tau|)^{-\og}d\tau\leq N,~~t\in\R;
$$

\item[\tu{(a$_3$)}] there exist positive constants $c$ and $\og$
such that, for any $\la,\la_1,\la_2\in Y$,
$$
\ba{l}
\|B(t,\la)\|\leq c
\min\{\mu(|t|)^{-\og-\ve},\nu(|t|)^{-\og-\ve}\},\\
\|B(t,\la_1)-B(t,\la_2)\|\leq c|\la_1-\la_2|\cdot\min\{\mu(|t|)^{-\og-\ve},\nu(|t|)^{-\og-\ve}\}.
\ea
$$
\end{itemize}
Moreover, if
\begin{equation}\label{eqthyangmi}
c<[KN(2K+1)]^{-1},
\end{equation}
then \eqref{eqpllinear} also admits a nonuniform
$(h,k,\mu,\nu)$-dichotomy on $\R$, i.e., for each $\la\in Y$, there exists a projection $\wpt(t,\la)$
for $t\in\R$ such that
\begin{equation}\label{invarianteo}
\wpt(t,\la)\htt(t,s,\la)=\htt(t,s,\la)\wpt(s,\la)
\end{equation}
and
\begin{equation}\label{boundes}
\begin{split}
&\|\htt(t,s,\la)\wpt(s,\la)\|\leq
\f{K\wkt}{1-2K\wkt cN}(h(t)/h(s))^a\mu(|s|)^{\ve}(\mu(|s|)^\ve+\nu(|s|)^\ve),~~t\geq s,\\
&\|\htt(t,s,\la)\wqt(s,\la)\|\leq \f{K\wkt}{1-2K\wkt
cN}(k(s)/k(t))^{-b}\nu(|s|)^{\ve}(\mu(|s|)^\ve+\nu(|s|)^\ve),~~s\geq
t,
\end{split}
\end{equation}
where $\wqt(t,\la)=\id-\wpt(t,\la)$ is the complementary projection of
$\wpt(t,\la)$ and $\htt(t,s,\la)$ is the evolution operator associated to \eqref{eqpllinear} and
\begin{equation}\label{eqthaaa}
\wkt=K/(1-KcN).
\end{equation}
Moreover, the stable
subspace $\wpt(t,\la)(X)$ and the unstable subspace
$\wqt(t,\la)(X)$ are Lipschitz continuous in $\la$
if $Y$ is finite-dimensional.
\end{theorem}

The proof of Theorem \ref{thaaa} is divided into several steps. First, we characterize the existence of some
bounded solutions of \eqref{eqpllinear} (Lemma \ref{leaaa1}) and
show that these bounded solutions are Lipschitz continuous in $\la$ (Lemma \ref{leaaa2}). We also show that these
bounded solutions admit a semigroup property for each $\la\in Y$
(Lemma \ref{lebbb2}). Then, with the help of the semigroup property
along bounded solutions and the invertible operator $S(0,\la)$
defined in Lemma \ref{leggg}, we construct invariant projections
$\wpt(t,\la)$ in \eqref{invarianteo}. Finally, we deliberately
establish the estimates for the evolution operator in
\eqref{boundes} (Lemmas \ref{lehhh}-\ref{lekkk}) and show that the
stable and unstable subspaces of nonuniform
$(h,k,\mu,\nu)$-dichotomies for the linear perturbed system are
Lipschitz continuous in $\la$ (Lemma \ref{lelll1}).

In the rest of this section, we always assume that the conditions in
Theorem \ref{thaaa} are satisfied.

{\bf Step 1.~Construction of bounded solutions of
\eqref{eqpllinear}.}

For each $s\in\R$, define
\begin{align*}
\Og_1&:=\{U(t,s)_{t\geq s}\in\B(X):U~\mbox{is continuous
and}~\|U\|_1<\iy,(t,s)\in\R\times\R\},\\
\Og_2&:=\{V(t,s)_{t\leq s}\in\B(X):V~\mbox{is continuous
and}~\|V\|_2<\iy,(t,s)\in\R\times\R\},
\end{align*}
respectively with the norms
\begin{align*}
&\|U\|_1=\sup\left\{\|U(t,s)\|(h(t)/h(s))^{-a}\mu(|s|)^{-\ve}:t\geq
s\right\},\\
&\|V\|_2=\sup\left\{\|V(t,s)\|(k(s)/k(t))^b\nu(|s|)^{-\ve}:t\leq
s\right\}.
\end{align*}
It is not difficult to show that $(\Og_1,\|\cdot\|_1)$ and
$(\Og_2,\|\cdot\|_2)$ are Banach spaces.

\begin{lemma}\label{leaaa1}
For each $\la\in Y$ and $s\in\R$,
\begin{itemize}
 \item there exists a unique solution $U^\la\in\Og_1$ of \eqref{eqpllinear} satisfying
\begin{equation}\label{eqzhubajie}
\begin{split}
U^\la(t,s)&=T(t,s)P(s)+\il_s^tT(t,\tau)P(\tau)B(\tau,\la)U^\la(\tau,s)d\tau\\
&\quad-\il_t^\iy
T(t,\tau)Q(\tau)B(\tau,\la)U^\la(\tau,s)d\tau,~~t\geq s;
\end{split}
\end{equation}

\item there exists a unique solution $V^\la\in\Og_2$ of
\eqref{eqpllinear} satisfying
\begin{equation}\label{eqwusong}
\begin{split}
V^\la(t,s)&=T(t,s)Q(s)+\il_{-\iy}^tT(t,\tau)P(\tau)B(\tau,\la)V^\la(\tau,s)d\tau\\
&\quad-\il_t^s
T(t,\tau)Q(\tau)B(\tau,\la)V^\la(\tau,s)d\tau,~~s\geq t.
\end{split}
\end{equation}

\end{itemize}
\end{lemma}
\begin{proof}
It is not difficult to show that $U^\la(t,s)_{t\geq s}$ satisfying
\eqref{eqzhubajie} and $V^\la(t,s)_{s\geq t}$ satisfying
\eqref{eqwusong} are solutions of \eqref{eqpllinear}. We next
show that the operator $J_1^\la$ defined by
\begin{align*}
(J_1^\la U)(t,s)&=T(t,s)P(s)+\il_s^tT(t,\tau)P(\tau)B(\tau,\la)U(\tau,s)d\tau\\
&\quad-\il_t^\iy T(t,\tau)Q(\tau)B(\tau,\la)U(\tau,s)d\tau
\end{align*}
has a unique fixed point in $\Og_1$ for each $\la\in Y$. For
$t\geq s$, by \eqref{deeqaaa}, (a$_2$) and (a$_3$), we have
\begin{align*}
A_1^\la:&=\il_s^t\|T(t,\tau)P(\tau)\|\|B(\tau,\la)\|\|U(\tau,s)\|d\tau+\il_t^\iy
\|T(t,\tau)Q(\tau)\|\|B(\tau,\la)\|\|U(\tau,s)\|d\tau\\
&\leq Kc(h(t)/h(s))^a\mu(|s|)^\ve\il_s^t\mu(|\tau|)^{-\og}
d\tau\|U\|_1\\
&\quad+Kc(h(t)/h(s))^{a}\mu(|s|)^\ve\il_t^\iy(k(\tau)/k(t))^{-b}(h(\tau)/h(t))^a\nu(|\tau|)^{-\og}
d\tau\|U\|_1\\
&\leq KcN(h(t)/h(s))^a\mu(|s|)^\ve\|U\|_1
\end{align*}
and
\begin{align*}
\|(J_1^\la U)(t,s)\|&\leq K(h(t)/h(s))^a\mu(|s|)^\ve+A_1^\la\\
&\leq
K(h(t)/h(s))^a\mu(|s|)^\ve+KcN(h(t)/h(s))^a\mu(|s|)^\ve\|U\|_1.
\end{align*}
Then
\begin{equation}\label{leeqfangshuan}
\|J_1^\la U\|_1\leq K+KcN\|U\|_1<\iy,
\end{equation}
which implies that $J_1^\la U$ is well-defined and
$J_1^\la:\Og_1\ra\Og_1$. Moreover, for each $\la\in Y$, for any
$U_1,U_2\in\Og_1$ and $t\geq s$, one has
\begin{align*}
A_2^\la:&=\il_s^t\|T(t,\tau)P(\tau)\|\|B(\tau,\la)\|\|U_1(\tau,s)-U_2(\tau,s)\|d\tau\\
&\leq Kc(h(t)/h(s))^a\mu(|s|)^\ve\il_s^t\mu(|\tau|)^{-\og}
d\tau\|U_1-U_2\|_1
\end{align*}
and
\begin{align*}
A_3^\la:&=\il_t^\iy\|T(t,\tau)Q(\tau)\|\|B(\tau,\la)\|\|U_1(\tau,s)-U_2(\tau,s)\|d\tau\\
&\leq
Kc(h(t)/h(s))^a\mu(|s|)^\ve\il_t^\iy\nu(|\tau|)^{-\og}d\tau\|U_1-U_2\|_1.
\end{align*}
Hence,
\begin{align*}
\|(J_1^\la U_1)(t,s)-(J_1^\la U_2)(t,s)\|\leq A_2^\la+A_3^\la\leq
KcN(h(t)/h(s))^a\mu(|s|)^\ve\|U_1-U_2\|_1,
\end{align*}
whence
$$
\|J_1^\la U_1-J_1^\la U_2\|_1\leq KcN\|U_1-U_2\|_1.
$$
If \eqref{eqthyangmi} holds, then the operator $J_1^\la$ is a
contraction and there exists a unique $U^\la\in\Og_1$ such that
$J_1^\la U^\la =U^\la.$

In addition, define an
operator $J_2^\la$ on $\Og_2$ by
\begin{align*}
(J_2^\la V)(t,s)&=T(t,s)Q(s)+\il_{-\iy}^tT(t,\tau)P(\tau)B(\tau,\la)V(\tau,s)d\tau\\
&\quad-\il_t^s T(t,\tau)Q(\tau)B(\tau,\la)V(\tau,s)d\tau
\end{align*}
for each $\la\in Y$. It follows from \eqref{deeqaaa}, (a$_2$), and (a$_3$) that
\begin{equation}\label{eqsongjiang}
\begin{split}
A_4^\la:&=\il_{-\iy}^t\|T(t,\tau)P(\tau)\|\|B(\tau,\la)\|\|V(\tau,s)\|d\tau+\il_t^s
\|T(t,\tau)Q(\tau)\|\|B(\tau,\la)\|\|V(\tau,s)\|d\tau\\
&\leq
Kc(k(s)/k(t))^{-b}\nu(|s|)^\ve\il_{-\iy}^t(h(t)/h(\tau))^{a}(k(t)/k(\tau))^{-b}\mu(|\tau|)^{-\og}d\tau\|V\|_2\\
&\quad+Kc(k(s)/k(t))^{-b}\nu(|s|)^\ve\il_t^s\nu(|\tau|)^{-\og}d\tau\|V\|_2\\
&\leq KcN(k(s)/k(t))^{-b}\nu(|s|)^\ve\|V\|_2
\end{split}
\end{equation}
and
\begin{align*}
\|(J_2^\la V)(t,s)\|&\leq K(k(s)/k(t))^{-b}\nu(|s|)^\ve+A_4^\la\\
&\leq
K(k(s)/k(t))^{-b}\nu(|s|)^\ve+KcN(k(s)/k(t))^{-b}\nu(|s|)^\ve\|V\|_2.
\end{align*}
Then
\begin{equation}\label{leeqfangshuanb}
\|J_2^\la V\|_2\leq K+KcN\|V\|_2<\iy
\end{equation}
and $J_2^\la:\Og_2\ra\Og_2$ is well-defined. On the other hand, we
have
$$
\|J_2^\la V_1-J_2^\la V_2\|_2\leq KcN\|V_1-V_2\|_2
$$
for each $\la\in Y$ and any $V_1,V_2\in\Og_2$. The operator $J_2^\la$ is a contraction since \eqref{eqthyangmi} holds
and there exists a unique $V^\la\in\Og_2$ such that $J_2^\la V^\la=V^\la$. The proof is complete.
\end{proof}

\begin{lemma}\label{leaaa2}
Both $U^\la$ and $V^\la$ are Lipschitz continuous in the
parameter $\la$.
\end{lemma}
\begin{proof}
It follows from Lemma \ref{leaaa1} that, for any $\la_1,\la_2\in Y$,
there exist bounded solutions $U^{\la_1},U^{\la_2}\in\Og_1$ satisfying \eqref{eqzhubajie}. Then, by \eqref{leeqfangshuan} and (a$_3$),
\begin{align*}
A^{\la_1,\la_2}_1(\tau):&=\|B(\tau,\la_1)U^{\la_1}(\tau,s)-B(\tau,\la_2)U^{\la_2}(\tau,s)\|\\
&\leq\|B(\tau,\la_1)U^{\la_1}(\tau,s)-B(\tau,\la_1)U^{\la_2}(\tau,s)\|+
\|B(\tau,\la_1)U^{\la_2}(\tau,s)-B(\tau,\la_2)U^{\la_2}(\tau,s)\|\\
&\leq
c(h(\tau)/h(s))^a\mu(|\tau|)^{-\og-\ve}\mu(|s|)^\ve(\|U^{\la_1}-U^{\la_2}\|_1+\wkt
|\la_1-\la_2|),
\end{align*}
which, together with \eqref{eqzhubajie}, implies
\begin{align*}
&\|U^{\la_1}(t,s)-U^{\la_2}(t,s)\|\\&\leq
\il_s^t\|T(t,\tau)P(\tau)\|A^{\la_1,\la_2}_1(\tau)d\tau+\il_t^\iy
\|T(t,\tau)Q(\tau)\|A^{\la_1,\la_2}_1(\tau)d\tau\\
&\leq
Kc(h(t)/h(s))^a\mu(|s|)^\ve\left(\il_s^t\mu(|\tau|)^{-\og}d\tau
+\il_t^\iy\nu(|\tau|)^{-\og}d\tau\right)
(\|U^{\la_1}-U^{\la_2}\|_1+\wkt |\la_1-\la_2|)\\
 &\leq KcN(h(t)/h(s))^a\mu(|s|)^\ve(\|U^{\la_1}-U^{\la_2}\|_1+\wkt
 |\la_1-\la_2|).
\end{align*}
Thus
$$
\|U^{\la_1}-U^{\la_2}\|_1\leq[\wkt
KcN/(1-KcN)]\cdot|\la_1-\la_2|.
$$
Similarly, for any $\la_1,\la_2\in Y$, there
exist bounded solutions $V^{\la_1}, V^{\la_2}\in\Og_2$ satisfying \eqref{eqwusong} and
\begin{align*}
A^{\la_1,\la_2}_2(\tau):&=\|B(\tau,\la_1)V^{\la_1}(\tau,s)-B(\tau,\la_2)V^{\la_2}(\tau,s)\|\\
&\leq\|B(\tau,\la_1)V^{\la_1}(\tau,s)-B(\tau,\la_1)V^{\la_2}(\tau,s)\|+
\|B(\tau,\la_1)V^{\la_2}(\tau,s)-B(\tau,\la_2)V^{\la_2}(\tau,s)\|\\
&\leq
c(k(\tau)/k(s))^{-b}\nu(|\tau|)^{-\og-\ve}\nu(|s|)^\ve(\|V^{\la_1}-V^{\la_2}\|_2+\wkt
|\la_1-\la_2|).
\end{align*}
Then
\begin{align*}
\|V^{\la_1}(t,s)-V^{\la_2}(t,s)\|&\leq
\il_{-\iy}^t\|T(t,\tau)P(\tau)\|A^{\la_1,\la_2}_2(\tau)d\tau+\il_t^s
\|T(t,\tau)Q(\tau)\|A^{\la_1,\la_2}_2(\tau)d\tau\\
&\leq KcN(k(t)/k(s))^a\nu(|s|)^\ve(\|V^{\la_1}-V^{\la_2}\|_2+\wkt
 |\la_1-\la_2|).
\end{align*}
The proof is complete.
\end{proof}

\emph{\bf Step 2.~Semigroup property of the bounded solutions.}

\begin{lemma}\label{lebbb2}
For each $\la\in Y$, one has
$$
U^\la(t,\sigma)U^\la(\sg,s)=U^\la(t,s),~t\geq\sg\geq s;~~V^\la(t,\sg)V^\la(\sg,s)=V^\la(t,s),~t\leq \sg\leq s.
$$
\end{lemma}
\begin{proof}
It follows from \eqref{eqzhubajie} that
\begin{align*}
U^\la(t,\sg)U^\la(\sg,s)&=T(t,s)P(s)+\il_s^\sg T(t,\tau)P(\tau)B(\tau,\la)U^\la(\tau,s)d\tau\\
&\quad+\il_\sg^tT(t,\tau)P(\tau)B(\tau,\la)U^\la(\tau,\sg)d\tau
U^\la(\sg,s)\\&\quad-\il_t^\iy T(t,\tau)Q(\tau)B(\tau,\la)U^\la(\tau,\sg)d\tau
U^\la(\sg,s).
\end{align*}
Let $L^\la(t,\sg)=U^\la(t,\sg)U^\la(\sg,s)-U^\la(t,s)$ for $t\geq \sg\geq s$.
Define the operator $H_1^\la$ by
\begin{align*}
(H_1^\la l)(t,\sg)=\il_\sg^tT(t,\tau)P(\tau)B(\tau,\la)l(\tau,\sg)d\tau-\il_t^\iy
T(t,\tau)Q(\tau)B(\tau,\la)l(\tau,\sg)d\tau,~~l\in\Og_1^\sg,~t\geq \sg,
\end{align*}
where $\Og_1^\sg$ is
obtained from $\Og_1$ by replacing $s$ with $\sg$. For any
$l,l_1,l_2\in\Og_1^\sg$, by (a$_1$)-(a$_3$), one has
\[
\|(H_1^\la l)(t,\sg)\|\leq KcN(h(t)/h(s))^a\mu(|s|)^\ve\|l\|_1
\]
and
\[
\|(H_1^\la l_1)(t,\sg)-(H_1^\la l_2)(t,\sg)\|\leq
KcN(h(t)/h(s))^a\mu(|s|)^\ve\|l_1-l_2\|_1,
\]
then
$$\|H_1^\la l\|_1\leq KcN\|l\|_1<\iy,~~\|H_1^\la l_1-H_1^\la l_2\|_1\leq KcN\|l_1-l_2\|_1.$$
Hence,
$H_1^\la$ is well-defined and $H_1^\la(\Og_1^\sg)\subset \Og_1^\sg$. Therefore, there
exists a unique $l^\la\in\Og_1^\sg$ such that $H_1^\la l^\la=l^\la$. Moreover, it
is not difficult to show that $L^\la\in\Og_1^\sg$ and $0\in\Og_1^\sg$ satisfying $H_1^\la0=0$ and $H_1^\la L^\la=L^\la$ , which implies that  $L^\la=l^\la=0$.
Similarly, by \eqref{eqwusong}, (a$_1$), (a$_2$), and (a$_3$), one has
$V^\la(t,\sg)V^\la(\sg,s)=V^\la(t,s)$ for $t\leq \sg\leq s$.
\end{proof}

\emph{\bf Step 3.~Construction of the projection $\wpt(t,\la)$ in
\eqref{invarianteo}.}

For each given $\la\in Y$. Define the linear operator
$$
\pt(t,\la)=\htt(t,0,\la)U^\la(0,0)\htt(0,t,\la),~~\qtt(t,\la)=\htt(t,0,\la)V^\la(0,0)\htt(0,t,\la),~t\in\R.
$$
By  Lemma \ref{lebbb2}, $\pt(t,\la)$ and $\qtt(t,\la)$ are projections for
each $t\in\R$ and
$$
\pt(t,\la)\htt(t,s,\la)=\htt(t,s,\la)\pt(s,\la),~~\qtt(t,\la)\htt(t,s,\la)=\htt(t,s,\la)\qtt(s,\la),~~t,s\in\R.
$$
It is obvious that $\widehat{U}^\la(t,0)=U^\la(t,0)P(0)$ satisfies \eqref{eqzhubajie} with $s=0$ and
$\widehat{V}^\la(t,0)=V^\la(t,0)Q(0)$ satisfies \eqref{eqwusong} with $s=0$. By Lemma \ref{leaaa1},
$$
U^\la(t,0)P(0)=U^\la(t,0),~~V^\la(t,0)Q(0)=V^\la(t,0).
$$
Note that
\begin{equation}\label{eqeqbaiyunruia}
\pt(0,\la)=U^\la(0,0)=P(0)-\il_0^\iy
T(0,\tau)Q(\tau)B(\tau,\la)U^\la(\tau,0)d\tau
\end{equation}
and
\begin{equation}\label{eqeqbaiyunruib}
\qtt(0,\la)=V^\la(0,0)=Q(0)+\il_{-\iy}^0
T(0,\tau)P(\tau)B(\tau,\la)V^\la(\tau,0)d\tau,
\end{equation}
then
\begin{equation}\label{step3formular}
\ba{c}
P(0)\pt(0,\la)=P(0),~\pt(0,\la)P(0)=\pt(0,\la),~P(0)(\id-\qtt(0,\la))=\id-\qtt(0,\la);\\
Q(0)\qtt(0,\la)=Q(0),~\qtt(0,\la)Q(0)=\qtt(0,\la),~Q(0)(\id-\pt(0,\la))=\id-\pt(0,\la).
\ea
\end{equation}
To obtain the projection $\widehat{P}(t,\la)$, set
$S(0,\la)=\pt(0,\la)+\qtt(0,\la)$.

\begin{lemma}\label{leggg}
For each $\la\in Y$, the operator $S(0,\la)$ is invertible.
\end{lemma}
\begin{proof}
It follows from \eqref{step3formular} that
\begin{equation}\label{leeqxuliang}
\pt(0,\la)+\qtt(0,\la)-\id=Q(0)\pt(0,\la)+P(0)\qtt(0,\la).
\end{equation}
By \eqref{eqeqbaiyunruia} and \eqref{eqeqbaiyunruib},
\begin{align*}
&P(0)\qtt(0,\la)=P(0)V^\la(0,0)=\il_{-\iy}^0
T(0,\tau)P(\tau)B(\tau,\la)V^\la(\tau,0)d\tau,\\
&Q(0)\pt(0,\la)=Q(0)U^\la(0,0)=-\il_0^\iy
T(0,\tau)Q(\tau)B(\tau,\la)U^\la(\tau,0)d\tau.
\end{align*}
Moreover, by \eqref{eqthaaa}, \eqref{leeqfangshuan} and
\eqref{leeqfangshuanb}, one has
\begin{equation}\label{leeqzhanzhaob}
\ba{l}
\|U^\la(t,s)\|\leq \wkt(h(t)/h(s))^a\mu(|s|)^\ve,~ t\geq s,\\
\|V^\la(t,s)\|\leq \wkt(k(s)/k(t))^{-b}\nu(|s|)^\ve, ~t\leq s.
\ea
\end{equation}
By \eqref{leeqxuliang}-\eqref{leeqzhanzhaob}, we have
\begin{align*}
A_5^\la&=:\il_{-\iy}^0
\|T(0,\tau)P(\tau)\|\|B(\tau,\la)\|\|V^\la(\tau,0)\|d\tau\\
&\leq K\wkt
c\il_{-\iy}^0(h(0)/h(\tau))^a(k(0)/k(\tau))^{-b}\mu(|\tau|)^{-\og}\nu(0)^\ve d\tau\\
&\leq K\wkt c\il_{-\iy}^0\mu(|\tau|)^{-\og}d\tau
\end{align*}
and
\begin{align*}
A_6^\la&=:\il_0^\iy
\|T(0,\tau)Q(\tau)\|\|B(\tau,\la)\|\|U^\la(\tau,0)\|d\tau\\
&\leq K\wkt
c\il_0^\iy(k(\tau)/k(0))^{-b}(h(\tau)/h(0))^a\mu(0)^\ve\nu(|\tau|)^{-\og}
d\tau\\
&\leq K\wkt
c\il_0^\iy\nu(|\tau|)^{-\og}d\tau.
\end{align*}
Then
\begin{align*}
\|\pt(0,\la)+\qtt(0,\la)-\id\|\leq A_5^\la+A_6^\la \leq K\wkt cN.
\end{align*}
Therefore, for each $\la\in Y$, the operator $S(0,\la)$ is invertible if
\eqref{eqthyangmi} holds.
\end{proof}

For $\la\in Y$ and $t\in\R$, set
\begin{equation}\label{eqliwenting}
\begin{split}
&\wpt(t,\la)=\htt(t,0,\la)S(0,\la)P(0)S(0,\la)^{-1}\htt(0,t,\la),\\
&\wqt(t,\la)=\htt(t,0,\la)S(0,\la)Q(0)S(0,\la)^{-1}\htt(0,t,\la).
\end{split}
\end{equation}
Then $\wpt(t,\la)$, $\wqt(t,\la)$ are projections for $t\in\R$ and $\wpt(t,\la)+\wqt(t,\la)=\id$.
Hence, \eqref{invarianteo} is valid.

\emph{\bf Step 4.~Norm bounds for the evolution operator.}

\begin{lemma}\label{lehhh} For each $\la\in Y$,
$\|\htt(t,s,\la)|\Imm\pt(s,\la)\|\leq \wkt(h(t)/h(s))^a\mu(|s|)^\ve$ for
$t\geq s$.
\end{lemma}
\begin{proof}
First, we show that, if $z^\la(t)_{(t\geq s)}$,  $\la\in Y$ is a bounded solution
of \eqref{eqpllinear} , then
\begin{equation}\label{eqsuiguangyu}
\begin{split}
z^\la(t)&=T(t,s)P(s)z^\la(s)+\il_s^tT(t,\tau)P(\tau)B(\tau,\la)z^\la(\tau)d\tau\\
&\quad-\il_t^\iy
T(t,\tau)Q(\tau)B(\tau,\la)z^\la(\tau)d\tau,~~t\geq s.
\end{split}
\end{equation}
It is not difficult to show that
\begin{equation}\label{eqyujiajia}
\ba{l}
P(t)z^\la(t)=T(t,s)P(s)z^\la(s)+\il_s^tT(t,\tau)P(\tau)B(\tau,\la)z^\la(\tau)d\tau,\\
Q(t)z^\la(t)=T(t,s)Q(s)z^\la(s)+\il_s^tT(t,\tau)Q(\tau)B(\tau,\la)z^\la(\tau)d\tau,
\ea
\end{equation}
and $z^\la(t)=P(t)z^\la(t)+Q(t)z^\la(t)$ for $t\in\R$. Then
\begin{equation}\label{eqbili}
Q(s)z^\la(s)=T(s,t)Q(t)z^\la(t)-\il_s^tT(s,\tau)Q(\tau)B(\tau,\la)z^\la(\tau)d\tau.
\end{equation}
On the other hand, one has
$$
\|T(s,t)Q(t)\|\leq K(k(t)/k(s))^{-b}\nu(|t|)^\ve
$$
and
\begin{align*}
\il_s^\iy\|T(s,\tau)Q(\tau)B(\tau,\la)z^\la(\tau)\|d\tau&\leq
Kc\il_s^\iy\nu(|\tau|)^{-\og}d\tau\sup\limits_{\tau\geq
s}\|z^\la(\tau)\|\\
&\leq KcN\sup\limits_{\tau\geq s}\|z^\la(\tau)\|<\iy.
\end{align*}
Let $t\ra \iy$ in \eqref{eqbili}, then
$$
Q(s)z^\la(s)=-\il_s^\iy T(s,\tau)Q(\tau)B(\tau,\la)z^\la(\tau)d\tau.
$$
Consequently,
\begin{align*}
Q(t)z^\la(t)&=-\il_s^\iy
T(t,\tau)Q(\tau)B(\tau,\la)z^\la(\tau)d\tau+\il_s^tT(t,\tau)Q(\tau)B(\tau,\la)z^\la(\tau)d\tau\\
&=-\il_t^\iy T(t,\tau)Q(\tau)B(\tau,\la)z^\la(\tau)d\tau,
\end{align*}
which proves \eqref{eqsuiguangyu}.

For each given $\xi\in X$ and  $\la\in Y$, let $z^\la(t)=\htt(t,s,\la)\pt(s,\la)\xi$ be the
solution of \eqref{eqpllinear} for $t\geq s$.
Since $\htt(t,0,\la)U^\la(0,0)$ and $U^\la(t,0)$ are solutions of
\eqref{eqpllinear} and coincide at $t=0$, then
$$
z^\la(t)=\htt(t,0,\la)U^\la(0,0)\htt(0,s,\la)\xi=U^\la(t,0)\htt(0,s,\la)\xi.
$$
Note that
$U^\la(t,0)$ is bounded for $t\in\R$, then
 $z^\la(t)_{(t\geq s)}$
is a bounded solution of \eqref{eqpllinear} with the initial value
$z^\la(s)=\pt(s,\la)\xi$. By \eqref{eqsuiguangyu},  we have
\begin{align*}
\pt(t,\la)\htt(t,s,\la)\xi&=T(t,s)P(s)\pt(s,\la)\xi+\il_s^tT(t,\tau)P(\tau)B(\tau,\la)\pt(\tau,\la)\htt(\tau,s,\la)\xi
d\tau\\
&\quad-\il_t^\iy T(t,\tau)Q(\tau)B(\tau,\la)\pt(\tau,\la)\htt(\tau,s,\la)\xi
d\tau,~~t\geq s.
\end{align*}
It is not difficult to show that
\begin{align*}
A_7^\la&=:\il_s^t\|T(t,\tau)P(\tau)\|\|B(\tau,\la)\|\|\pt(\tau,\la)\htt(\tau,s,\la)\xi\|
d\tau\\
&\leq
Kc\il_s^t(h(t)/h(\tau))^{a}\mu(|\tau|)^{-\og}\|\pt(\tau,\la)\htt(\tau,s,\la)\|\|\pt(s,\la)\xi\|d\tau\\
&\leq
Kc(h(t)/h(s))^a\mu(|s|)^\ve\|\pt(\la)\htt(\la)\|_1\|\pt(s,\la)\xi\|\il_s^t\mu(|\tau|)^{-\og}d\tau
\end{align*}
and
\begin{align*}
A_8^\la&=:\il_t^\iy
\|T(t,\tau)Q(\tau)\|\|B(\tau,\la)\|\|\pt(\tau,\la)\htt(\tau,s,\la)\xi
\|d\tau\\
&\leq
Kc\il_t^\iy(k(\tau)/k(t))^{-b}\nu(|\tau|)^{-\og}\|\pt(\tau,\la)\htt(\tau,s,\la)\|\|\pt(s,\la)\xi\|d\tau\\
&\leq Kc\|\pt(\la)\htt(\la)\|_1\il_t^\iy(k(\tau)/k(t))^{-b}\nu(|\tau|)^{-\og}(h(\tau)/h(s))^a\mu(|s|)^\ve\|\pt(s,\la)\xi\|d\tau\\
&\leq Kc(h(t)/h(s))^a\mu(|s|)^\ve\|\pt(\la)\htt(\la)\|_1\|\pt(s,\la)\xi\|\il_t^\iy\nu(|\tau|)^{-\og}d\tau.
\end{align*}
Then
\begin{align*}
\|\pt(t,\la)\htt(t,s,\la)\xi\|&\leq
K(h(t)/h(s))^a\mu(|s|)^\ve\|\pt(s,\la)\xi\|+A_7^\la+A_8^\la\\
&\leq K(h(t)/h(s))^a\mu(|s|)^\ve\|\pt(s,\la)\xi\|\\
&\quad+
KcN(h(t)/h(s))^a\mu(|s|)^\ve\|\pt(\la)\htt(\la)\|_1\|\pt(s,\la)\xi\|,
\end{align*}
i.e., $\|\pt(\la)\htt(\la)\|_1\leq \wkt.$ This yields the desired
inequality.
\end{proof}

\begin{lemma}\label{leiii}
For each $\la\in Y$, $\|\htt(t,s,\la)|\Imm\qtt(s,\la)\|\leq\wkt(k(s)/k(t))^{-b}\nu(|s|)^\ve$ for
$t\leq s$.
\end{lemma}

\begin{proof}
By carrying out arguments similar to that of Lemma \ref{lehhh}, we can show that, for $\la\in Y$, if $z^\la(t)_{(t\leq s)}$ is a bounded
solution of \eqref{eqpllinear}, then
\begin{equation}\label{eqwanghongyan}
\begin{split}
z^\la(t)&=T(t,s)Q(s)z^\la(s)+\il_{-\iy}^tT(t,\tau)P(\tau)B(\tau,\la)z^\la(\tau)d\tau\\
&\quad-\il_t^sT(t,\tau)Q(\tau)B(\tau,\la)z^\la(\tau)d\tau.
\end{split}
\end{equation}
Moreover,
\begin{align*}
z^\la(t):=\htt(t,s,\la)\qtt(s,\la)\xi=V^\la(t,0)\htt(0,s,\la)\xi, ~~\xi\in X,~t\leq s
\end{align*}
and $z^\la(t)_{(t\leq s)}$ is a bounded
solution of \eqref{eqpllinear} with $z^\la(s)=\qtt(s,\la)\xi$. From \eqref{eqwanghongyan}, it follows that
\begin{align*}
\qtt(t,\la)\htt(t,s,\la)\xi&=T(t,s)Q(s)\qtt(s,\la)\xi+\il_{-\iy}^tT(t,\tau)P(\tau)B(\tau,\la)\qtt(\tau,\la)\htt(\tau,s,\la)\xi
d\tau\\
&\quad-\il_t^s T(t,\tau)Q(\tau)B(\tau,\la)\qtt(\tau,\la)\htt(\tau,s,\la)\xi
d\tau.
\end{align*}
Note that
\begin{align*}
A_9^\la&=:\il_{-\iy}^t\|T(t,\tau)P(\tau)\|\|B(\tau,\la)\|\|\qtt(\tau,\la)\htt(\tau,s,\la)\xi\|
d\tau\\
&\leq Kc\il_{-\iy}^t(h(t)/h(\tau))^{a}\mu(\tau)^{-\og}
\|\qtt(\tau,\la)\htt(\tau,s,\la)\|\|\qtt(s,\la)\xi\|d\tau\\
&\leq
Kc\|\qtt(\la)\htt(\la)\|_2\il_{-\iy}^t(h(t)/h(\tau))^{a}\mu(|\tau|)^{-\og}(k(s)/k(\tau))^{-b}\nu(|s|)^\ve
\|\qtt(s,\la)\xi\|d\tau\\
&\leq
Kc(k(s)/k(t))^{-b}\nu(|s|)^\ve\|\qtt(\la)\htt(\la)\|_2\|\qtt(s,\la)\xi\|\il^t_{-\iy}\mu(|\tau|)^{-\og}d\tau
\end{align*}
and
\begin{align*}
A_{10}^\la&=:\il_t^s\|T(t,\tau)Q(\tau)\|\|B(\tau,\la)\|\|\qtt(\tau,\la)\htt(\tau,s,\la)\xi\|
d\tau\\
&\leq Kc\il_t^s(k(\tau)/k(t))^{-b}\nu(|\tau|)^{-\og}
\|\qtt(\tau,\la)\htt(\tau,s,\la)\|\|\qtt(s,\la)\xi\|d\tau\\
&\leq Kc(k(s)/k(t))^{-b}\nu(|s|)^\ve\|\qtt(\la)\htt(\la)\|_2\|\qtt(s,\la)\xi\|\il_t^s\nu(|\tau|)^{-\og}d\tau,
\end{align*}
we have
\begin{align*}
\|\qtt(t,\la)\htt(t,s,\la)\xi\|&\leq
K(k(s)/k(t))^{-b}\nu(|s|)^\ve\|\qtt(s,\la)\xi\|+A_9^\la+A_{10}^\la\\
&\leq
K(k(s)/k(t))^{-b}\nu(|s|)^\ve\|\qtt(s,\la)\xi\|\\
&\quad+KcN(k(s)/k(t))^{-b}\nu(|s|)^\ve\|\qtt(\la)\htt(\la)\|_2\|\qtt(s,\la)\xi\|.
\end{align*}
Then
$$
\|\qtt(\la)\htt(\la)\|_2\leq K +KcN\|\qtt(\la)\htt(\la)\|_2,~~\hbox{i.e.,}~~
\|\qtt(\la)\htt(\la)\|_2\leq \wkt,
$$ which yields the desired inequality.
\end{proof}

\begin{lemma}\label{lejjj}
For each $\la\in Y$, one has
\begin{equation}\label{eqooppqq}
\begin{split}
&\|\htt(t,s,\la)\wpt(s,\la)\|\leq\wkt(h(t)/h(s))^a\mu(|s|)^\ve\|\wpt(s,\la)\|,~~t\geq
s,\\
&\|\htt(t,s,\la)\wqt(s,\la)\|\leq\wkt(k(s)/k(t))^{-b}\nu(|s|)^\ve\|\wqt(s,\la)\|,~~t\leq
s.
\end{split}
\end{equation}
\end{lemma}

\begin{proof}
For $\la\in Y$, by (b$_4$), we have
\begin{align*}
&S(0,\la)P(0)=(\pt(0,\la)+\qtt(0,\la))P(0)=\pt(0,\la),\\&S(0,\la)Q(0)=(\pt(0,\la)+\qtt(0,\la))Q(0)=\qtt(0,\la)
\end{align*}
Note that $S(t,\la)=\htt(t,0,\la)S(0,\la)\htt(0,t,\la)$ for $t\in\R$, then
\begin{align*}
\wpt(t,\la)S(t,\la)=\htt(t,0,\la)S(0,\la)P(0)\htt(0,t,\la)=\htt(t,0,\la)\pt(0,\la)\htt(0,t,\la)=\pt(t,\la).
\end{align*}
Similarly, $\wqt(t,\la)S(t,\la)=\qtt(t,\la)$. Then
\begin{align*}
\Imm\wpt(t,\la)=\Imm\pt(t,\la)\tx{and}\Imm\wqt(t,\la)=\Imm\qtt(t,\la).
\end{align*}
By Lemmas \ref{lehhh} and \ref{leiii}, one has
\begin{align*}
\|\htt(t,s,\la)\wpt(s,\la)\|&\leq\|\htt(t,s,\la)|\Imm\wpt(s,\la)\|\|\wpt(s,\la)\|\\
&\leq\|\htt(t,s,\la)|\Imm\pt(s,\la)\|\|\wpt(s,\la)\|\\
&\leq\wkt(h(t)/h(s))^a\mu(|s|)^\ve\|\wpt(s,\la)\|,~~t\geq s
\end{align*}
and
\begin{align*}
\|\htt(t,s,\la)\wqt(s,\la)\|&\leq\|\htt(t,s,\la)|\Imm\wqt(s,\la)\|\|\wqt(s,\la)\|\\
&\leq\|\htt(t,s,\la)|\Imm\qtt(s,\la)\|\|\wqt(s,\la)\|\\
&\leq\wkt(k(s)/k(t))^{-b}\nu(|s|)^\ve\|\wqt(s,\la)\|,~~t\leq s.
\end{align*}
\end{proof}

\begin{lemma}\label{lekkk}
For each $\la\in Y$, one has
\begin{equation}\label{eqrroott}
\begin{split}
\|\wpt(t,\la)\|&\leq
[K/(1-2K\wkt cN)](\mu(|t|)^\ve+\nu(|t|)^\ve),\\
\|\wqt(t,\la)\|&\leq [K/(1-2K\wkt cN)](\mu(|t|)^\ve+\nu(|t|)^\ve).
\end{split}
\end{equation}
\end{lemma}

\begin{proof}
For $\xi\in X$ and $\la\in Y$, set
\begin{align*}
z_1^\la(t)=\htt(t,s,\la)\wpt(s,\la)\xi,~t\geq
s;~~z_2^\la(t)=\htt(t,s,\la)\wqt(s,\la)\xi,~t\leq s.
\end{align*}
By Lemma \ref{lejjj}, $(z_1^\la(t))_{t\geq s}$ and $(z_2^\la(t))_{t\leq
s}$ are bounded solutions of \eqref{eqpllinear}. By
\eqref{eqsuiguangyu} and \eqref{eqwanghongyan},
\begin{align*}
\wpt(t,\la)\htt(t,s,\la)\xi&=T(t,s)P(s)\wpt(s,\la)\xi+\il_s^tT(t,\tau)P(\tau)B(\tau,\la)\wpt(\tau,\la)\htt(\tau,s,\la)\xi
d\tau\\
&\quad-\il_t^\iy T(t,\tau)Q(\tau)B(\tau,\la)\wpt(\tau,\la)\htt(\tau,s,\la)\xi
d\tau
\end{align*}
and
\begin{align*}
\wqt(t,\la)\htt(t,s,\la)\xi&=T(t,s)Q(s)\wqt(s,\la)\xi+\il_{-\iy}^tT(t,\tau)P(\tau)B(\tau,\la)\wqt(\tau,\la)\htt(\tau,s,\la)\xi
d\tau\\
&\quad-\il_t^s T(t,\tau)Q(\tau)B(\tau,\la)\wqt(\tau,\la)\htt(\tau,s,\la)\xi
d\tau.
\end{align*}
Taking $t=s$ leads to
\begin{align*}
&Q(t)\wpt(t,\la)\xi=-\il_t^\iy
T(t,\tau)Q(\tau)B(\tau,\la)\wpt(\tau,\la)\htt(\tau,t,\la)\xi d\tau,\\
&P(t)\wqt(t,\la)\xi=\il_{-\iy}^tT(t,\tau)P(\tau)B(\tau,\la)\wqt(\tau,\la)\htt(\tau,t,\la)\xi
d\tau.
\end{align*}
By Lemma \ref{lejjj},
\begin{align*}
\|Q(t)\wpt(t,\la)\|+\|P(t)\wqt(t,\la)\| &\leq \wkt Kc\left(\mu(|t|)^\ve\il_t^\iy\nu(|\tau|)^{-\og}d\tau\|\wpt(t,\la)\|\right.\\
&\quad\left.
+\nu(|t|)^\ve\il_{-\iy}^t\mu(|\tau|)^{-\og}d\tau\|\wqt(t,\la)\|\right)\\
&\leq\wkt KcN(\|\wpt(t,\la)\|+\|\wqt(t,\la)\|).
\end{align*}
Since $\|P(t)\|\leq K \mu(|t|)^\ve$ and $\|Q(t)\|\leq K \nu(|t|)^\ve$,
one has
\begin{align*}
\|\wpt(t,\la)\|&\leq\|\wpt(t,\la)-P(t)\|+\|P(t)\|=\|\wpt(t,\la)-P(t)\wpt(t,\la)-P(t)+P(t)\wpt(t,\la)\|+\|P(t)\|\\
&=\|Q(t)\wpt(t,\la)-P(t)\wqt(t,\la)\|+\|P(t)\|\leq\|Q(t)\wpt(t,\la)\|+\|P(t)\wqt(t,\la)\|+\|P(t)\|\\
&\leq\wkt KcN(\|\wpt(t,\la)\| +\|\wqt(t,\la)\|)+K\mu(|t|)^\ve
\end{align*}
and
\begin{align*}
\|\wqt(t,\la)\|&\leq\|\wqt(t,\la)-Q(t)\|+\|Q(t)\|=\|\wpt(t,\la)-P(t)\|+\|Q(t)\|\\
&\leq \wkt KcN(\|\wpt(t,\la)\| +\|\wqt(t,\la)\|)+K\nu(|t|)^\ve.
\end{align*}
Therefore,
$$
\|\wpt(t,\la)\|+\|\wqt(t,\la)\|\leq 2\wkt KcN (\|\wpt(t,\la)\|
+\|\wqt(t,\la)\|)+K(\mu(|t|)^\ve+\nu(|t|)^\ve).
$$
The proof is complete.
\end{proof}

\emph{\bf Step 5.~Lipschitz continuity of $\wpt(t,\la)(X)$,
$\wqt(t,\la)(X)$ with respect to $\la$.
}
\begin{lemma}\label{lelll1}
$\wpt(t,\la)(X)$ and
$\wqt(t,\la)(X)$ are Lipschitz continuous in $\la$.
\end{lemma}
\begin{proof}
By (a$_3$), $\htt(t,0,\la)$ is Lipschitz continuous in $\la$.Since
$U^{\lambda}$ and $V^{\lambda}$ are Lipschitz continuous in $\lambda$ (Lemma \ref{leaaa2}),
  $\pt(t,\lambda)$ and $\qtt(t,\lambda)$ are Lipschitz continuous in $\lambda$. Moreover, if $Y$ is
finite-dimensional, then $S(0,\la)$ and $S^{-1}(0,\la)$ are
both Lipschitz continuous in the parameter. By \eqref{eqliwenting}, Lemma \ref{lelll1} is valid.
\end{proof}

\begin{remark}
Theorem \ref{thxwq} includes and generalizes Theorem 4.1 in \cite{Barreira2013}
(nonuniform $(\mu,\nu)$-dichotomies) and Theorem 7 in
\cite{Barreira2009bb} (nonuniform exponential dichotomy). When \eqref{eqpllinear} reduces to $x'=(A(t)+B(t))x$, the
conclusion in Theorem \ref{thaaa} includes and extends some
existing results for robustness of various dichotomies, such as, robustness of exponential
dichotomy (Theorem 5.6 in \cite{plh2006}, Theorem 3.2 in
\cite{Ju2001} and Proposition 1 of Section 4 in
\cite{Coppelbook1978}), robustness of  $(h,k)$-dichotomy (Theorem
6 in \cite{Naulin1995}), robustness of  nonuniform exponential
dichotomy (Theorem 2 in \cite{Barreira200855}), robustness of
$\rho$-nonuniform exponential dichotomy (Theorem 2 in
\cite{Barreira2009b}), and robustness of nonuniform
$(\mu,\nu)$-dichotomy (Theorem 4.1 in \cite{Chang2011}).
\end{remark}

\section{Existence of topological conjugacy}\label{ghsection}
\noindent

In this section, with the help of the nonuniform
$(h,k,\mu,\nu)$-dichotomy, we explore the topological conjugacy of
nonautonomous dynamical systems in Banach spaces by establishing a
new version of the Grobman-Hartman theorem.

Consider the
nonlinear perturbed system
\begin{equation}\label{eqnonlinear}
x'=A(t)x+f(t,x),
\end{equation}
where $f:\R\times X\ra X$.
\begin{definition}[see \cite{Palmer1973}]\label{defintionanas;ldf} \eqref{eqlinear} and \eqref{eqnonlinear} are said to
be topologically equivalent if there exists an function
$H:\R\times X\ra X$ having the following properties:
\begin{itemize}
\item[\tu{(i)}] if $\|x\|\ra\iy$, then $\|H(t,x)\|\ra\iy$
uniformly with respect to $t\in\R$;

\item[\tu{(ii)}] for each fixed $t$, $H(t,\cdot)$ is a
homeomorphism of $X$ into $X$;

\item[\tu{(iii)}] $L(t,\cdot)=H^{-1}(t,\cdot)$ also has property
(i);

\item[\tu{(iv)}] if $x(t)$ is a solution of \eqref{eqnonlinear},
then $H(t,x(t))$ is a solution of \eqref{eqlinear}.
\end{itemize}

The function $H$ satisfying the above four properties is said to
be \emph{the equivalent function} of \eqref{eqlinear} and
\eqref{eqnonlinear}.
\end{definition}

\begin{theorem}\label{ghtheorem}
Assume that \eqref{eqlinear} admits a nonuniform
$(h,k,\mu,\nu)$-dichotomy as in \eqref{deeqaaa} on $\R$ and $h,k$
are  differentiable. If there exist positive constants $\al$ and
$\ga$ such that, for any $x,x_1,x_2\in X$,
\begin{equation}\label{eqliuxiu}
\begin{split}
\|f(t,x)\|&\leq\al\min\{h'(t)h(t)^{-1}\mu(|t|)^{-\ve},k'(t)k(t)^{-1}\nu(|t|)^{-\ve}\},\\
\|f(t,x_1)-f(t,x_2)\|&\leq
\ga\min\{h'(t)h(t)^{-1}\mu(|t|)^{-\ve},k'(t)k(t)^{-1}\nu(|t|)^{-\ve}\}\|x_1-x_2\|,
\end{split}
\end{equation}
 and
\begin{equation}\label{eqwangmang}
K\ga(1/|a|+1/b)<1,
\end{equation}
then \eqref{eqnonlinear} is topologically equivalent to
\eqref{eqlinear} and the equivalent function $H(t,x)$ satisfies
$$
\|H(t,x)-x\|\leq K\al(1/|a|+1/b),~~t\in\R,~~x\in X.
$$
\end{theorem}

The proof of Theorem \ref{ghtheorem} is achieved in two steps.
First, it is shown that either of the systems \eqref{eqzmzaa},
\eqref{eqzwaa} and \eqref{eqswk} has a unique bounded solution
(Lemma \ref{lezbj}, \ref{leldh}, \ref{lesss}), then we construct a
function $H(t,x)$ and prove that $H(t,x)$ is an equivalent
function satisfying the properties (i)-(iv) in Definition
\ref{defintionanas;ldf} (Lemma \ref{leyjj}, \ref{lejcy},
\ref{leabc}, \ref{lexyz}).

Let
$X(t,t_0, x_0)$ be the solution of \eqref{eqnonlinear} with
$X(t_0)=x_0$ and $Y(t,t_0,y_0)$ be the solution of
\eqref{eqlinear} with $Y(t_0)=y_0$. In the rest of this section, we always assume that \eqref{eqliuxiu} and \eqref{eqwangmang} are satisfied.

{\bf Step 1.~Construction of bounded solutions.}

\begin{lemma}\label{lezbj}
For any fixed $(\bar{t}, \xi)\in\R\times X$,
\begin{equation}\label{eqzmzaa}
z'=A(t)z-f(t,X(t,\bar{t},\xi))
\end{equation}
has a unique bounded solution $h(t,(\bar{t},\xi))$ and
$$\|h(t,(\bar{t},\xi))\|\leq
K\al(1/|a|+1/b),~~t\in\R.$$
\end{lemma}

\begin{proof}
It is trivial to show that
\begin{align*}
h(t,(\bar{t},\xi))&=-\il_{-\iy}^tT(t,\tau)P(\tau)f(\tau,X(\tau,\bar{t},\xi))d\tau
+\il_t^\iy T(t,\tau)Q(\tau)f(\tau,X(\tau,\bar{t},\xi))d\tau
\end{align*}
is a solution of \eqref{eqzmzaa}. By \eqref{deeqaaa} and
\eqref{eqliuxiu}, for any $t\in\R$, we have
\begin{align*}
\|h(t,(\bar{t},\xi))\|&=\il_{-\iy}^t\|T(t,\tau)P(\tau)\|\|f(\tau,X(\tau,\bar{t},\xi))\|d\tau\\
&\quad+\il_t^\iy \|T(t,\tau)Q(\tau)\|\|f(\tau,X(\tau,\bar{t},\xi))\|d\tau\\
&\leq K\al h(t)^a\il_{-\iy}^th(\tau)^{-a-1}h'(\tau)d\tau+K\al h(t)^b\il_t^\iy h(\tau)^{-b-1}h'(\tau)d\tau\\
&\leq K\al(1/|a|+1/b),
\end{align*}
which implies that $h(t,(\bar{t},\xi))$ is the unique bounded
solution of \eqref{eqzmzaa} since $z'=A(t)z$ admits a nonuniform
$(h,k,\mu,\nu)$-dichotomy on $\R$.
\end{proof}

\begin{lemma}\label{leldh}
For any fixed $(\bar{t},\xi)\in\R\times X$,
\begin{equation}\label{eqzwaa}
z'=A(t)z+f(t,Y(t,\bar{t},\xi)+z)
\end{equation}
has a unique bounded solution $l(t,(\bar{t},\xi))$ and
$$\|l(t,(\bar{t},\xi))\|\leq
K\al(1/|a|+1/b).$$
\end{lemma}

\begin{proof}
Let
$$
\Og_3:=\{z:\R\ra X|\|z\|\leq K\al(1/|a|+1/b)\},
$$
where $\|z\|:=\sup_{t\in\R}\|z(t)\|$. Then
 $(\Og_3,\|\cdot\|)$ is a Banach
space. Define a mapping  $J$ on $\Og_3$ by
\begin{align*}
(Jz)(t)&=\il_{-\iy}^tT(t,\tau)P(\tau)f(\tau,Y(\tau,\bar{t},\xi)+z(\tau))d\tau\\
&\quad-\il_t^\iy
T(t,\tau)Q(\tau)f(\tau,Y(\tau,\bar{t},\xi)+z(\tau))d\tau.
\end{align*}
It follows from \eqref{deeqaaa} and \eqref{eqliuxiu} that
$$
\|Jz\|\leq K\al(1/|a|+1/b),~~\|Jz_1-Jz_2\|\leq
K\ga(1/|a|+1/b)\|z_1-z_2\|,~~z,z_1,z_2\in\Og_3.
$$
Then $J(\Og_3)\subset\Og_3$ and $J$ is a contraction
mapping. Therefore, $J$ has a unique fixed point
$l(t)$, i.e.,
\begin{align*}
l(t,(\bar{t},\xi))&=\il_{-\iy}^tT(t,\tau)P(\tau)f(\tau,Y(\tau,\bar{t},\xi)+l(\tau))d\tau\\
&\quad-\il_t^\iy
T(t,\tau)Q(\tau)f(\tau,Y(\tau,\bar{t},\xi)+l(\tau))d\tau.
\end{align*}
Next, we prove that $l(t,(\bar{t},\xi))$ is unique in
the whole space by contradiction arguments. Otherwise, assume that there is another bounded
solution $l^0(t,(\bar{t},\xi))$ of \eqref{eqzwaa}, which can be
written as
\begin{align*}
l^0(t,(\bar{t},\xi))&=\il_{-\iy}^tT(t,\tau)P(\tau)f(\tau,Y(\tau,\bar{t},\xi)+l^0(\tau))d\tau\\
&\quad-\il_t^\iy
T(t,\tau)Q(\tau)f(\tau,Y(\tau,\bar{t},\xi)+l^0(\tau))d\tau.
\end{align*}
It is trivial to show that
$$
\|l-l^0\|\leq K\ga(1/|a|+1/b)\|l-l^0\|.
$$
Then, by \eqref{eqwangmang}, one has $l\equiv l^0$. Therefore,
$l(t,(\bar{t},\xi))$ is a unique bounded solution of
\eqref{eqzwaa} and
$$\|l(t,(\bar{t},\xi))\|\leq
K\al(1/|a|+1/b),~t\in\R.$$
\end{proof}

\begin{lemma}\label{lesss}
Let $x(t)$ be any solution of \eqref{eqnonlinear}, then
\begin{equation}\label{eqswk}
z'=A(t)z+f(t,x(t)+z)-f(t,x(t))
\end{equation}
has a unique bounded solution $z(t)\equiv0$.
\end{lemma}
\begin{proof}
It is obvious that $z(t)\equiv0$ is a bounded solution of
\eqref{eqswk}. Next we show that $z(t)\equiv0$ is the unique
bounded solution. Assume that $z^0(t)$ is any bounded solution of
\eqref{eqswk}, then $z^0(t)$ can be written in the form
\begin{align*}
z^0(t)&=\il_{-\iy}^t T(t,\tau)P(\tau)[f(\tau,x(\tau)+z(\tau))-f(\tau,x(\tau))]d \tau\\
&\quad-\il_t^\iy
T(t,\tau)Q(\tau)[f(\tau,x(\tau)+z(\tau))-f(\tau,x(\tau))]d\tau.
\end{align*}
It is easy to show that
$$
\|z^0-0\|\leq K\ga(1/|a|+1/b)\|z^0-0\|,
$$
which implies that $z^0(t)\equiv0$.
\end{proof}

{\bf Step 2.~Construction of the topologically equivalent
function.}

Define
\begin{equation}\label{eqwly}
H(t,x)=x+h(t,(t,x)),~~L(t,y)=y+l(t,(t,y)),~~x,y\in X.
\end{equation}
\begin{lemma}\label{leyjj}
For any fixed $(\bar{t},x(\bar{t}))\in\R\times X$,
$H(t,X(t,\bar{t},x(\bar{t})))$ is a solution of \eqref{eqlinear}.
\end{lemma}
\begin{proof}
By Lemma \ref{lezbj}, we have
$$
h(t,(t,X(t,\bar{t},x(\bar{t}))))=h(t,(\bar{t},x(\bar{t})))
$$
and
$$
H(t,X(t,\bar{t},x(\bar{t})))=X(t,\bar{t},x(\bar{t}))+h(t,(t,X(t,\bar{t},x(\bar{t}))))
=X(t,\bar{t},x(\bar{t}))+h(t,(\bar{t},x(\bar{t}))).
$$
Note that $X(t,\bar{t},x(\bar{t}))$ and
$h(t,(\bar{t},x(\bar{t})))$ are solutions of \eqref{eqnonlinear}
and \eqref{eqzmzaa}, respectively, then
\begin{align*}
H'(t,X(t,\bar{t},x(\bar{t})))&=X'(t,\bar{t},x(\bar{t}))+h'(t,(\bar{t},x(\bar{t})))\\
&=A(t)X(t,\bar{t},x(\bar{t}))+f(t,X(t,\bar{t},x(\bar{t})))\\
&\quad+A(t)h(t,(\bar{t},x(\bar{t})))-f(t,X(t,\bar{t},x(\bar{t})))\\
&=A(t)H(t,X(t,\bar{t},x(\bar{t}))),
\end{align*}
which implies that  $H(t,X(t,\bar{t},x(\bar{t})))$ is a solution of
\eqref{eqlinear}.
\end{proof}

\begin{lemma}\label{lejcy}
For any fixed $(\bar{t}, y(\bar{t}))\in\R\times X$,
$L(t,Y(t,\bar{t},y(\bar{t})))$ is a solution of
\eqref{eqnonlinear}.
\end{lemma}
\begin{proof}
It follows from Lemma \ref{leldh} that
$$
l(t,(t,Y(t,\bar{t},y(\bar{t}))))=l(t,(\bar{t},y(\bar{t}))),
$$
then
\begin{align*}
L(t,Y(t,\bar{t},y(\bar{t})))&=Y(t,\bar{t},y(\bar{t}))+l(t,(t,Y(t,\bar{t},y(\bar{t}))))\\
&=Y(t,\bar{t},y(\bar{t}))+l(t,(\bar{t},y(\bar{t}))).
\end{align*}
Since $Y(t,\bar{t},y(\bar{t}))$ and $l(t,(\bar{t},y(\bar{t})))$
are solutions of \eqref{eqlinear} and \eqref{eqzwaa},
respectively,  we have
\begin{align*}
L'(t,Y(t,\bar{t},y(\bar{t})))&=Y'(t,\bar{t},y(\bar{t}))+l'(t,(\bar{t},y(\bar{t})))\\
&=A(t)Y(t,\bar{t},y(\bar{t}))+A(t)l(t,(\bar{t},y(\bar{t})))\\&\quad+f(t,Y(t,\bar{t},y(\bar{t}))+l(t,(\bar{t},y(\bar{t}))))\\
&=A(t)L(t,Y(t,\bar{t},y(\bar{t})))+f(t,L(t,Y(t,\bar{t},y(\bar{t})))).
\end{align*}
\end{proof}

\begin{lemma}\label{leabc}
For any fixed $t\in\R$ and $y\in X$, $H(t,L(t,y))=y$ holds.
\end{lemma}

\begin{proof}
Let $y(t)$ be any solution of \eqref{eqlinear}. It follows from
Lemma \ref{leyjj} and
 Lemma \ref{lejcy} that $L(t,y(t))$ is a
solution of \eqref{eqnonlinear} and $H(t,L(t,y(t)))$ is a solution
of \eqref{eqlinear}. Moreover,
\begin{align*}
H'(t,L(t,y(t)))-y'(t)&=A(t)H(t,L(t,y(t)))-A(t)y(t)=A(t)(H(t,L(t,y(t)))-y(t))
\end{align*}
and
\begin{align*}
\|H(t,L(t,y(t)))-y(t)\|&\leq\|H(t,L(t,y(t)))-L(t,y(t))\|+\|L(t,y(t))-y(t)\|\\
&\leq2K\al(1/|a|+1/b).
\end{align*}
Then $H(t,L(t,y(t)))-y(t)$ is a bounded solution of
\eqref{eqlinear} and
$H(t,L(t,y(t)))-y(t)\equiv0.$
For any fixed $t\in\R, y\in X$, there is a solution of
\eqref{eqlinear} with the initial value $y(t)=y$. Then
$H(t,L(t,y))=y$ holds.
\end{proof}

\begin{lemma}\label{lexyz}
For any fixed $t\in\R$ and  $x\in X$, $L(t,H(t,x))=x$ holds.
\end{lemma}

\begin{proof}
Let $x(t)$ be any solution of \eqref{eqnonlinear}. It follows from
Lemma \ref{leyjj} and Lemma \ref{lejcy} that $H(t,x(t))$ is a
solution of \eqref{eqlinear} and $L(t,H(t,x(t)))$ is a solution of
\eqref{eqnonlinear}. Moreover,
$$
\ba{l}
L'(t,H(t,x(t)))-x'(t)\\
~~~~~~~~=A(t)L(t,H(t,x(t)))+f(t,L(t,H(t,x(t))))-A(t)x(t)-f(t,x(t))\\
~~~~~~~~=A(t)[L(t,H(t,x(t)))-x(t)]+f(t,L(t,H(t,x(t)))-x(t)+x(t))-f(t,x(t))
\ea
$$
and
\begin{align*}
\|L(t,H(t,x(t)))-x(t)\|&\leq\|L(t,H(t,x(t)))-H(t,x(t))\|+\|H(t,x(t))-x(t)\|\\
&\leq2K\al(1/|a|+1/b).
\end{align*}
By Lemma \ref{lesss}, we conclude that
$L(t,H(t,x(t)))-x(t)\equiv0.$ For any fixed $t\in\R, x\in X$,
there exists a solution of \eqref{eqnonlinear} with the initial
value $x(t)=x$. Then $L(t,H(t,x))=x$ holds.
\end{proof}

We are now at the right position to establish Theorem
\ref{ghtheorem}, that is, to verify that $H(t,x)$ is topologically
equivalent function. From \eqref{eqwly} and Lemma \ref{lezbj}, it
follows that, for any $t\in\R$,
$$\|H(t,x)-x\|=\|h(t,x)\|\leq
K\al(1/|a|+1/b),~ x\in X.$$ Then $\|H(t,x)\|\ra\iy$ as
$\|x\|\ra\iy$ uniformly with respect to $t\in\R$, i.e., Condition
(i) holds. By Lemma \ref{leabc} and Lemma \ref{lexyz}, for each
fixed $t\in\R$, $H(t,\cdot)=L^{-1}(t,\cdot)$ is homeomorphism.
Then Condition (ii) holds. By \eqref{eqwly} and Lemma \ref{leldh},
for any $t\in\R$, we have
$$\|L(t,y)-y\|=\|l(t,,y)\|\leq
K\al(1/|a|+1/b),~y\in X.$$ This implies that $\|L(t,y)\|\ra\iy$ as
$\|y\|\ra\iy$ uniformly with respect to $t\in\R$. Hence, Condition
(iii) holds. It follows from Lemma \ref{leyjj} and Lemma
\ref{lejcy} that Condition (iv) holds.

\begin{remark}
Theorem \ref{ghtheorem} not only includes the
classical Palmer's linearization  theorem for hyperbolic system in
\cite{Palmer1979}, and also extends the idea of linearization theorems
from hyperbolicity to nonuniform hyperbolicity.
\end{remark}

\section{Existence of stable invariant manifolds}\label{manifoldsection}
\noindent

We establish in this section the existence of parameter dependence of Lipschitz stable
invariant manifolds for sufficiently small Lipschitz perturbations
of \eqref{eqlinear} assuming that it admits a nonuniform
$(h,k,\mu,\nu)$-dichotomy.

Consider the nonlinear perturbed system with the parameters of \eqref{eqlinear}
\begin{equation}\label{eqnonlinearaaa}
x'=A(t)x+f(t,x,\la),
\end{equation}
where $f:\R\times X\times Y\ra X$ and $f(t,0,\la)=0$ for any $t\in\R$ and $\la\in Y$. Since the problem explored here is the existence of stable invariant manifold of \eqref{eqnonlinearaaa},
one only needs to carry out the discussion on $\R^+$.
In order to facilitate the discussion below, we make use of the following equivalent characterization of the nonuniform $(h,k,\mu,\nu)$-dichotomy
\begin{equation}\label{deeqaaabbb}
\begin{split}
&\|T(t,s)P(s)\|\le K\left(\f{h(t)}{h(s)}\right)^a\mu(s)^\ve,~~\|T(t,s)^{-1}Q(t)\| \le
K\left(\f{k(t)}{k(s)}\right)^{-b}\nu(t)^\ve,~~t\geq s\geq 0.
\end{split}
\end{equation}
Assume that there exist positive constants
$\hc$ and $q$ such that
\begin{equation}\label{*A1}
\begin{split}
\|f(t,x_1,\la)-f(t,x_2,\la)\|&\le
\hat{c}\|x_1-x_2\|(\|x_1\|^{q}+\|x_2\|^{q}),\\
\|f(t,x,\la_1)-f(t,x,\la_2)\|&\le
\hat{c}|\la_1-\la_2|\cdot\|x\|^{q+1}
\end{split}
\end{equation}
for any $t\in\R^+$,$x, x_1, x_2\in X$ and $\la,\la_1,\la_2\in Y$. Define the stable and
unstable spaces for each $t\in\R^+$ by
$$
E(t) =P(t) (X) \quad \text{and} \quad F(t) =Q(t) (X).
$$
We next establish the parameter dependence of the stable manifolds as graphs of Lipschitz
functions and begin with introducing the class of functions to be
considered. For each $s\ge 0$, let $B_s(\varrho)\subset E(s)$ be
the open ball of radius $\varrho$ centered at zero and set
\begin{equation}\label{*beta}
\bt(t)=k(t)^{b/(\ve q)}h(t)^{-a(q+1)/(\ve
q)}\mu(t)^{1+1/q}C(t)^{1/\ve q},
\end{equation}
where
$$
C(t)=\il_t^\iy h(\tau)^{aq}\max\{\mu(\tau)^\ve,\nu(\tau)^\ve\}
d\tau.
$$
Given $\eta>0$, consider the set of initial conditions
\[
Z_\bt(\eta)=\big\{(s,\xi):s\geq0,~ \xi\in
B_s(\bt(s)^{-\ve}/\eta)\big\}.
\]
Let $Z_\bt =Z_\bt(1)$. Denote by $\Xa$ the space of continuous functions
$\Phi\colon Z_\bt\to X$ such that
\[
\Phi(s,0)=0, \quad \Phi(s,B_s(\bt(s)^{-\ve}))\subset F(s),
\]
and
\begin{equation}\label{*cuca}
\|\Phi(s,\xi_1)-\Phi(s,\xi_2)\|\le\|\xi_1-\xi_2\|
\end{equation}
for $s \ge 0$ and $\xi_1,\xi_2\in B_s(\bt(s)^{-\ve})$. It is
not difficult to show that $\Xa$ is a complete metric space
induced by
\[
|\Phi|'=\sup\left\{\frac{\|\Phi(s,\xi)\|}{\|\xi\|}: s\geq0
~\mbox{and}~ \xi\in B_s(\bt(s)^{-\ve})\setminus \{0\}\right\}.
\]
For each $\la\in Y$ and given $\Phi \in \Xa$, consider the graph
\begin{equation}\label{*ru}
\W^\la=\big\{(s,\xi,\Phi(s,\xi)): (s,\xi)\in Z_\bt\big\}
\end{equation}
and the semiflow generated by \eqref{eqnonlinearaaa}:
\begin{equation}\label{*psi}
\Psi_\kappa^\la(s,u(s), v(s))=(t, u(t), v(t)),~~\kappa=t-s\geq 0,~(s,u(s),v(s))\in\R^+\times E(s)\times F(s)
\end{equation}
where
\begin{equation}\label{eqbl}
\ba{l}
u(t)=T(t,s)u(s)+
\il_s^{t}T(t,\tau)P(\tau)f(\tau,u(\tau),v(\tau),\la)d\tau,\\
v(t)=T(t,s)v(s)+\il_s^{t}T(t,\tau)Q(\tau)f(\tau,u(\tau),v(\tau),\la)d\tau.
\ea
\end{equation}

We now state the existence of parameter dependence of a stable invariant manifold for
\eqref{eqnonlinearaaa}.

\begin{theorem}\label{theoremlsm}
Assume that

\begin{itemize}
\item[\tu{(c$_1$)}] \eqref{eqlinear} admits a nonuniform
$(h,k,\mu,\nu)$-dichotomy as in \eqref{deeqaaa} on $\R^+$;

\item[\tu{(c$_2$})]
$\lim\limits_{t\ra\iy}k(t)^{-b}h(t)^a\nu(t)^\ve=0;$

\item[\tu{(c$_3$)}] $h(t)^a\bt(t)^\ve$ is decreasing.
\end{itemize}
If  $\hat{c}$ in \eqref{*A1} is sufficiently small, then for each $\la\in Y$,
\begin{itemize}
\item[\tu{(d$_1$)}] there exists a unique function $\Phi=\Phi^\la\in \Xa$
such that $\W^\la$ is forward invariant with respect to $\Psi_\kal^\la$ in
the sense that
\begin{equation}\label{eqkkk}
\Psi_\kappa^\la(s,\xi,\Phi(s,\xi))\in \W^\la \quad \text{for any} \quad
(s,\xi)\in Z_{\bt\cdot\mu}(2K), ~\kappa=t-s\geq0;
\end{equation}

\item[\tu{(d$_2$)}]there exists a constant $d>0$ such that
\begin{equation}\label{*sec}
\|\Psi_\kappa^\la(s,\xi_1,\Phi(s,\xi_1))-\Psi_\kappa^\la(s,\xi_2,\Phi(s,\xi_2))\|
\le d(h(t)/h(s))^a\mu(s)^\ve\|\xi_1-\xi_2\|
\end{equation}
for any $\kappa=t-s \ge 0 $ and $(s,\xi_1),(s,\xi_2) \in
Z_{\bt\cdot\mu}(2K)$;

\item[\tu{(d$_3$)}] there exists a constant $d^*>0$ such that
\begin{equation}\label{*secb}
\|\Psi_\kappa^{\la_1}(t,\xi,\Phi^{\la_1}(t,\xi))-\Psi_\kappa^{\la_2}(t,\xi,\Phi^{\la_2}(t,\xi))\|
\le d(h(t)/h(s))^a\mu(s)^\ve|\la_1-\la_2|\cdot\|\xi\|.
\end{equation}
for any $\la_1,\la_2\in Y$.
\end{itemize}
\end{theorem}

To obtain parameter dependence of the stable manifolds, we first introduce an auxiliary
space. Let $\bar{\Xa}$ be the space of functions
$\Phi\colon\R^+\times X\to X$ such that $\Phi|_{Z_\bt}\in X$ and
\[
\Phi(s,\xi)=\Phi\big(s, \bt(s)^{-\ve}\xi/\|\xi\|\big),~~(s,\xi) \not \in Z_\bt.
 \]
Note that there is a one-to-one correspondence
between $\Xa$ and $\bar{\Xa}$. Moreover,
$\bar{\Xa}$ is a Banach space with the norm
$\bar{\Xa}\ni\Phi\mapsto|\Phi|Z_\bt|'$. It is not difficult to
show that, for each $\Phi\in\bar{\Xa}$, one has
\begin{equation}\label{eqzmz}
\|\Phi(s,\xi_1)-\Phi(s,\xi_2)\|\le2\|\xi_1-\xi_2\|,~~s \ge 0,~\xi_1, \xi_2 \in E(s).
\end{equation}

The proof of Theorem \ref{theoremlsm} is obtained in several
steps. First, we prove that, for each  $(s,\xi,\Phi,\la)\in Z_\bt\times\Xa^*\times Y$, there
exists a unique function $u^{\Phi,\la}(t,s,\xi)$ satisfying the first equality of \eqref{eqbl} (Lemma
\ref{leaaa}). In order to prove that there exists a unique
function $\Phi=\Phi^\la\in \Xa$ satisfying \eqref{eqzwt} for each $\la\in Y$, we reduce the problem to an
alternative one (Lemma \ref{leaabb}) and show that there exists a
unique function $\Phi=\Phi^\la\in \Xa$ satisfying \eqref{eqyjj} (Lemma \ref{leopq}). The asymptotic behavior of the
unique function $u^{\Phi,\la}(t,s,\xi)$ defined in Lemma \ref{leaaa} are characterized by Lemma \ref{lezx}, \ref{lemma.4}, and \ref{le123456}. Finally, with the established lemmas, we prove Theorem
\ref{theoremlsm} by showing that \eqref{eqkkk}, \eqref{*sec}, and \eqref{*secb}
are satisfied.

\begin{lemma}\label{leaaa}
Let $\hat{c}$ in \eqref{*A1} be sufficiently small. Then, for each $(s,\xi,\Phi,\la)\in Z_\bt\times\Xa^*\times Y$, there exists a unique function
$u=u^{\Phi,\la}:\R^+\to X$ with $u(s)=\xi$ such that, for any $t\geq s$,
\eqref{eqbl} holds and
\begin{equation}\label{eqcxy}
\|u(t)\|\le 2 K (h(t)/h(s))^a\mu(s)^\ve\|\xi\|.
\end{equation}
\end{lemma}

\begin{proof}
Let $\Omega_4$ be the space of continuous functions $u \colon
[s,\iy)\to X$ with $u(s)=\xi$ such that $u(t)\in E(t)$ for $t \ge s$ and $\|u\|_*\le \bt(s)^{-\ve}$, where
\begin{equation}\label{eq:norma}
\|u\|_*=\frac{1}{2K}\sup\left\{\frac{\|u(t)\|}{(h(t)/h(s))^a\mu(s)^\ve}:
t\geq s\right\}.
\end{equation}
It is trivial to show that $\Og_4$ is a complete metric space
induced by $\lVert \cdot\rVert_*$.

Given $(s,\xi)\in Z_\bt$ and $\Phi\in\bar{\Xa}$, for $t \ge s$ and each $\la\in Y$,
define an operator $L^\la$ in $\Og_4$ by
\begin{align*}
(L^\la u)(t)&=T(t,s)\xi+\il_s^t T(t,\tau)P(\tau)f(\tau,u(\tau),
\Phi(\tau,u(\tau),\la)d\tau.
\end{align*}
Then $L^\la u$ is continuous in $[s,\iy)$,
$(L^\la u)(s)=\xi$ and $(L^\la u)(t)\in E(t)$ for $t \ge s$. It
follows from \eqref{deeqaaabbb} and  \eqref{*A1} that
\begin{align*}
B_1^\la(\tau)&=:\|f(\tau,u(\tau),\Phi(\tau,u(\tau)),\la)\|\\
&\le \hc\left(\|u(\tau)\|+\|\Phi(\tau,u(\tau))\|\right)
\left(\|u(\tau)\|+\|\Phi(\tau,u(\tau))\|\right)^{q}\\
&\le 3^{q+1}\hc\|u(\tau)\|^{q+1}\\&\le 6^{q+1}\hc
K^{q+1}\left(\f{h(\tau)}{h(s)}\right)^{
a(q+1)}\mu(s)^{\ve(q+1)}(\|u\|_*)^{q+1},~\tau\geq s
\end{align*}
and
\begin{align*}
\|(L^\la u)(t)\|&\le\|T(t,s)\|\|\xi\|+\il_s^t \|T(t,\tau)P(\tau)\| B_1^\la(\tau) d\tau\\
&\le  K\left(\f{\mu(t)}{\mu(s)}\right)^a\nu(s)^\ve\|\xi\|+6^{q+1}\hc
K^{q+2}\lf\f{h(t)}{h(s)}\rf^ah(s)^{-aq}\mu(s)^{\ve(q+1)}(\|u\|_*)^{q+1}C(s),
\end{align*}
which implies that
\begin{align*}
\|L^\la u\|_*&\leq \f{1}{2}\left(\|\xi\|+6^{q+1}\hat{c}K^{q+1}
h(s)^{-aq}\mu(s)^{\ve q}(\|u\|_*)^{q+1}C(s)\right)\\
&\leq \f{1}{2}\left(1+6^{q+1}\hat{c}K^{q+1} h(s)^{-aq}\mu(s)^{\ve
q}\bt(s)^{-\ve q}C(s)
\right)\bt(s)^{-\ve}\leq\f{1}{2}(1+6^{q+1}\hat{c}K^{q+1})\bt(s)^{-\ve}.
\end{align*}
Since $\hat{c}$ is sufficiently small, take $\hat{c}$ such that
$6^{q+1}\hat{c}K^{q+1}< 1$, then $L^\la(\Og_4)\subset\Og_4$.
In addition, for any $u_1,u_2 \in\Og_4$, one has
\begin{align*}
B_2^\la(\tau)&=:\|f(\tau,u_1(\tau),\Phi(\tau,u_1(\tau)),\la)
-f(\tau,u_2(\tau),\Phi(\tau,u_2(\tau)),\la)\|\\
&\leq
3^{q+1}\hc\|u_1(\tau)-u_2(\tau)\|(\|\mu_1(\tau)\|^q+\|\mu_2(\tau)\|^q)\\
&\leq2^{q+2}3^{q+1}\hc
K^{q+1}\lf\f{h(\tau)}{h(s)}\rf^{a(q+1)}\mu(s)^{\ve(q+1)}\bt(s)^{-\ve
q}\|u_1-u_2\|_*
\end{align*}
and
\begin{align*}
\|L^\la u_1(t)-L^\la u_2(t)\|&\leq
\il_s^t\|T(t,\tau)P(\tau)\| B_2^\la(\tau) d\tau\\
&\leq2\cdot6^{q+1}\hat{c}K^{q+1}\|u_1-u_2\|_*
\lf\f{h(t)}{h(s)}\rf^a\mu(s)^\ve.
\end{align*}
Whence,
$$
\|L^\la u_1-L^\la u_2\|_*\leq6^{q+1}\hat{c}K^{q+1}\|u_1-u_2\|_*.
$$
Therefore, $L^\la$ is a contraction in $\Og_4$ and there exists a unique
function $u=u^{\Phi,\la}\in\Omega_4$ such that $L^u=u$. Moveover,
\[
\|u\|_*\le\f{1}{2}\|\xi\|+\f{1}{2}6^{q+1}\hat{c}K^{q+1}\|u\|_*,
\]
and
$$
\|u(t)\|\le 2 K (h(t)/h(s))^a\mu(s)^\ve\|\xi\| \quad \text{for
any} \quad t\geq s,
$$
since $K/(1-(1/2)6^{q+1}\hat{c}K^{q+1})< 2 K$.
\end{proof}

Let $u(t)=u^{\Phi,\la}(t,s,\xi)$ be the unique function defined by Lemma
\ref{leaaa}, that is,
\begin{equation}\label{*coc}
\begin{split}
u(t)&= T(t,s)\xi+\il_s^t T(t,\tau)P(\tau)f(\tau,u(\tau),
\Phi(\tau,u(\tau),\la)d\tau,~~t\geq s.
\end{split}
\end{equation}

\begin{lemma}\label{leaabb}
Given $\hat{c}>0$ sufficiently small and $\Phi\in \bar{\Xa}$, for each $\la\in Y$, the
following properties hold:
\begin{itemize}
\item[\tu{(e$_1$)}] for each $(s,\xi)\in Z_\bt$ and $t\ge s$, if
\begin{equation}\label{eqzwt}
\Phi(t,u(t))=T(t,s)\Phi(s,\xi)
+\il_s^tT(t,\tau)Q(\tau)f(\tau,u(\tau),\Phi(\tau,u(\tau)),\la)d\tau,
\end{equation}
then
\begin{equation}\label{eqyjj}
\Phi(s,\xi)=-\il_s^\iy
T(\tau,s)^{-1}Q(\tau)f(\tau,u(\tau),\Phi(\tau,u(\tau)),\la)d\tau;
\end{equation}
\item[\tu{(e$_2$)}] if \eqref{eqyjj} holds for $s \ge 0$ and
$\xi\in B_s(\bt(s)^{-\ve})$, then \eqref{eqzwt} holds for
$(s,\xi)\in Z_{\bt\cdot\mu}(2K)$.
\end{itemize}
\end{lemma}

\begin{proof}
By \eqref{deeqaaabbb}, \eqref{*A1},  \eqref{eqzmz}, and
\eqref{eqcxy}, for $\tau\ge s$, one has
\begin{align*}
B_3^\la(\tau)&=:\|T(\tau,s)^{-1}Q(\tau)\|\cdot \|f(\tau,u(\tau),\Phi(\tau,u(\tau)),\la)\|\\
&\le 3^{q+1}\hc K\left(\f{k(\tau)}{k(s)}\right)^{-b}\nu(\tau)^\ve\|u(\tau)\|^{q+1}\\
&\le6^{q+1}\hc K^{q+2}\left(\f{k(\tau)}{k(s)}\right)^{-b}
\nu(\tau)^\ve\lf\f{h(\tau)}{h(s)}\rf^{a(q+1)}\mu(s)^{\ve(q+1)}\|\xi\|^{q+1}\\
&\le 6^{q+1}\hc K^{q+2}\left(\f{k(\tau)}{k(s)}\right)^{-b}
\nu(\tau)^\ve\lf\f{h(\tau)}{h(s)}\rf^{a(q+1)}\mu(s)^{\ve(q+1)}\bt(s)^{-\ve(q+1)}
\end{align*}
and
\begin{align*}
\il_s^\iy B_3^\la(\tau)d\tau &\le 6^{q+1}\hc
K^{q+2}k(s)^{b}h(s)^{-a(q+1)}\mu(s)^{\ve(q+1)}\bt(s)^{-\ve(q+1)}\times\lf\il_s^\iy k(\tau)^{-b}h(\tau)^{a(q+1)}\nu(\tau)^\ve
d\tau\rf\\
&\leq6^{q+1}\hc
K^{q+2}k(s)^{b}h(s)^{-a(q+1)}\mu(s)^{\ve(q+1)}\bt(s)^{-\ve
q}C(s)<\iy,
\end{align*}
which imply that the right-hand side of \eqref{eqyjj} is
well-defined.

Assume that \eqref{eqzwt} holds for $(s,\xi)\in Z_\bt$ and $t\ge
s$, then \eqref{eqzwt} rewrites
\begin{equation}\label{*co}
\begin{split}
\Phi(s,\xi)&=T(t,s)^{-1}\Phi(t,u(t))
-\il_s^tT(\tau,s)^{-1}Q(\tau)f(\tau,u(\tau),\Phi(\tau,u(\tau)),\la)d\tau.
\end{split}
\end{equation}
From \eqref{deeqaaabbb}, \eqref{eqzmz}, and \eqref{eqcxy}, it
follows that
\begin{align*}
\|T(t,s)^{-1}\Phi(t,u(t))\|&\le4K^2\left(\f{k(t)}{k(s)}\right)^{-b}\nu(t)^\ve
\left(\f{h(t)}{h(s)}\right)^a\mu(s)^\ve\bt(s)^{-\ve}\\
&\le4K^2k(t)^{-b}h(t)^a\nu(t)^\ve
k(s)^{b}h(s)^{-a}\mu(s)^\ve\bt(s)^{-\ve}.
\end{align*}
Then, letting $t\to\iy$ in \eqref{*co} yields \eqref{eqyjj}.

If \eqref{eqyjj} holds for any $(s,\xi)\in Z_\beta$, then, for
$(s,\xi) \in Z_{\bt\cdot\mu}(2 K)$,
\[
\|u(t)\|\le2K
\left(\f{h(t)}{h(s)}\right)^a\mu(s)^\ve\|\xi\|\leq\bt(t)^{-\ve}\f{h(t)^a\bt(t)^{\ve}}{h(s)^a\bt(s)^{\ve}}
\leq\bt(t)^{-\ve},
\]
and hence, $(t,u(t))\in Z_\bt$ for any $t\ge s$. By \eqref{eqyjj},
one has
\begin{align*}
T(t,s)\Phi(s,\xi)
&=-\il_s^tT(t,\tau)Q(\tau)f(\tau,u(\tau),\Phi(\tau,u(\tau)),\la)d\tau\\
&\quad -\il_t^\iy
T(t,\tau)Q(\tau)f(\tau,u(\tau),\Phi(\tau,u(\tau)),\la)d\tau
\\&=-\il_s^t
T(t,\tau)Q(\tau)f(\tau,u(\tau),\Phi(\tau,u(\tau)),\la)d\tau+\Phi(t,\mu(t)),
\end{align*}
where, in the last equality, \eqref{eqyjj} is used with $(s,\xi)$ replaced
by $(t,u(t))$.
\end{proof}

\begin{lemma}\label{lezx}
If $\hat{c}$ is sufficiently small and $u_i(t)=u^{\Phi,\la}(t,s,\xi_i), i=1,2$, then there exists a $K_1>0$
such that
\begin{equation}\label{eqzw}
\|u_1(t)-u_2(t)\|\le K_1(h(t)/h(s))^a\mu(s)^\ve\|\xi_1-\xi_2\|,~~t\geq s.
\end{equation}
\end{lemma}

\begin{proof}
It follows from \eqref{*A1},
\eqref{eqzmz} and \eqref{eqcxy} that
\begin{align*}
B_4^\la(\tau)&=:\|f(\tau,u_1(\tau),
\Phi(\tau,u_1(\tau)),\la)-f(\tau,u_2(\tau),\Phi(\tau,u_2(\tau)),\la)\|\\
&\leq3^{q+1}\hc\|u_1(\tau)-u_2(\tau)\|(\|u_1(\tau)\|^q+\|u_2(\tau)\|^q)
\end{align*}
and
\begin{align*}
\|u_1(t)-u_2(t)\|&\leq\|T(t,s)(\xi_1 -\xi_2)\|+\il_s^t\|T(t,\tau)P(\tau)\|B_4^\la(\tau)d\tau\\
&\leq K\lf\f{h(t)}{h(s)}\rf^a\mu(s)^\ve(\|\xi_1-\xi_2\|\\
&\quad+2\cdot6^{q+1}\hc K^{q+2}\|u_1-u_2\|_*\times \lf\f{h(t)}{h(s)}\rf^a\mu(s)^{\ve(q+1)}\bt(s)^{-\ve
q}\il_s^th(\tau)^{aq}\nu(\tau)^\ve d\tau.
\end{align*}
Then
\[
\|u_1-u_2\|^*\le\f{1}{2}\|\xi_1-\xi_2\|+6^{q+1}\hat{c}K^{q+1}\|u_1-u_2\|^*,
\]
which yields~\eqref{eqzw} with $K_1=K/(1-6^{q+1}\hat{c}K^{q+1}).$
\end{proof}

\begin{lemma}\label{lemma.4}
If $\hat{c}$ is sufficiently small and let $u_i(t)=u^{\Phi_i,\la}(t,s,\xi), i=1,2$, then there exists a $K_2>0$
such that
\begin{equation}\label{eqwww}
\|u_1(t)-u_2(t)\|\le
K_2(\mu(t)/\mu(s))^a\|\xi\|\cdot |\Phi_1-\Phi_2|',~~t\geq s.
\end{equation}
\end{lemma}

\begin{proof}
For simplicity, write $u_i=u_\xi^{\Phi_i}$ for $i=1,2$. A
straightforward calculation shows that
\begin{align*}
B_5^\la(\tau)&=:\|f(\tau,u_1(\tau),\Phi_1(\tau,u_1(\tau)),\la)
-f(\tau,u_2(\tau),\Phi_2(\tau,u_2(\tau)),\la)\|\\
&\leq3^{q}\hc\left[3(\|u_1(\tau)-u_2(\tau)\|)
(\|u_1(\tau)\|^{q}+\|u_2(\tau)\|^{q})\right.\\
&\quad~~~~~+\left.(\|u_1(\tau)\|\cdot |\Phi_1-\Phi_2|')
(\|u_1(\tau)\|^{q}+\|u_2(\tau)\|^{q})\right]\\
&\leq[2\cdot6^{q+1}\hc K^{q+1}\|u_1-u_2\|_*\\
&~~~~~~~~+4\cdot6^q\hc K^{q+1}\|\xi\|\cdot|\Phi_1-\Phi_2|']\times \left(\f{h(\tau)}{h(s)}\right)^{a(q+1)}\mu(s)^{\ve(q+1)}\bt(s)^{-\ve
q},~\tau\geq s
\end{align*}
and
\begin{align*}
\|u_1(t)-u_2(t)\|&\leq\il_s^t\|T(t,\tau)P(\tau)\|B_5^\la(\tau)d\tau\\
&\leq[2\cdot6^{q+1}\hc K^{q+1}\|u_1-u_2\|_*+4\cdot6^q\hc K^{q+1}\|\xi\|\cdot|\Phi_1-\Phi_2|']\\
&~~~~~~\times
K\lf\f{h(t)}{h(s)}\rf^ah(s)^{-aq}\mu(s)^{\ve(q+1)}\bt(s)^{-\ve
q}C(s).
\end{align*}
Then
\begin{align*}
\|u_1-u_2\|_*\leq[6^{q+1}\hc K^{q+1}\|u_1-u_2\|_*+2\cdot6^q\hc
K^{q+1}\|\xi\|\cdot|\Phi_1-\Phi_2|']\mu(s)^{-\ve}.
\end{align*}
This establishes \eqref{eqwww}.
\end{proof}

\begin{lemma}\label{leopq}
If $\hat{c}$ is sufficiently small, then, for each $\la\in Y$, there exists a unique
function $\Phi=\Phi^\la\in \Xa$ such that \eqref{eqyjj} holds for any
$(s,\xi)\in Z_\bt$.
\end{lemma}

\begin{proof}
For each $\la\in Y$ and any $(s,\xi)\in Z_\bt$, define an operator
$J^\la$ by
\begin{align*}
(J^\la\Phi)(s,\xi)&=-\il_s^\iy
T(\tau,s)^{-1}Q(\tau)f(\tau,u(\tau),\Phi(\tau,u(\tau)),\la)d\tau,~~\Phi\in\bar{\Xa}
\end{align*}
where $u=u^{\Phi,\la}(t,s,\xi)$ is the unique function defined by
Lemma~\ref{leaaa}. It is not difficult to show that $J^\la\Phi$ is
continuous and $(J^\la\Phi)(s,0)=0$ for $s\ge 0$. Moreover, for
any $\xi_1, \xi_2\in B_s(\bt(s)^{-\ve})$, let
$u_i(t)=u^{\Phi,\la}(t,s,\xi_i)$, $i=1,2$,. By \eqref{deeqaaabbb},
\eqref{eqcxy}, and \eqref{eqzw}, we have
\begin{align*}
B_6^\la(\tau):&=\|f(\tau,u_1(\tau),\Phi(\tau,u_1(\tau)),\la)
-f(\tau,u_2(\tau),\Phi(\tau,u_2(\tau)),\la)\|\\
&\leq3^{q+1} \hc\| u_1(\tau)-u_2(\tau) \| (\| u_1(\tau)\|^q + \| u_2(\tau)\|^q )\\
&\leq6^{q+1}\hc
K^qK_1\lf\f{h(\tau)}{h(s)}\rf^{a(q+1)}\mu(s)^{\ve(q+1)}\bt(s)^{-\ve
q}\|\xi_1-\xi_2\|
\end{align*}
and
\begin{align*}
\|(J^\la\Phi)(s,\xi_1)-(J^\la\Phi)(s,\xi_2)\|&\leq\il_s^\iy\|T(\tau,s)^{-1}Q(\tau)\|B_6^\la(\tau)d\tau\\
&\leq6^{q+1}\hc K^{q+1}K_1k(s)^{b}h(s)^{-a(q+1)}\times\mu(s)^{\ve(q+1)}\bt(s)^{-\ve q}C(s)\|\xi_1-\xi_2\|\\
&\leq6^{q+1}cK^{q+1}K_1\|\xi_1-\xi_2\|.
\end{align*}
If $\hat{c}$ is sufficiently small, then
\[
\|(J^\la\Phi)(s,\xi_1)-(J^\la\Phi)(s,\xi_2)\|\le\|\xi_1-\xi_2\|
\]
and one can extend $J \Phi$ to $\R^+ \times X$ by $
(J\Phi)(s,\xi)= (J\Phi)\big(s,\bt(s)^{-\ve}\xi /\lVert \xi \rVert
\big) $ for any $(s,\xi) \not \in Z_\bt$. Hence,
$J(\bar{\Xa})\subset \bar{\Xa}$.

Now we show that $J$ is a contraction. For any $\Phi_1,\Phi_2\in
\bar{\Xa}$, write $u_i(t)=u^{\Phi_i,\la}(t,s,\xi)$, $i=1,2$, for each
$(s,\xi)\in Z_\bt$, by
\eqref{eqzmz}, \eqref{eqcxy}, and \eqref{eqwww},  one has
\begin{align*}
B_7^\la(\tau):&=\|f(\tau,u_1(\tau), \Phi_1(\tau,u_1(\tau)),\la)-
f(\tau,u_2(\tau),\Phi_2(\tau,u_2(\tau)),\la)\|\\
&\leq3^q\hc\big(3\|u_1(\tau)-u_2(\tau)\|+
\|u_1(\tau)\|\cdot|\Phi_1-\Phi_2|'\big)
(\|u_1(\tau)\|^q+\|u_2(\tau)\|^q)\\
&\leq2\cdot6^q\hc K^{q}(2K+3K_2)\|\xi\|\cdot |\Phi_1-\Phi_2|'\times\lf\f{h(\tau)}{h(s)}\rf^{a(q+1)}\mu(s)^{\ve(q+1)}\bt(s)^{-\ve
q}
\end{align*}
and
\begin{align*}
\|(J^\la\Phi_1)(s,\xi)-(J^\la\Phi_2)(s,\xi)\|
&\leq\il_s^\iy\|T(\tau,s)^{-1}Q(\tau)\|B_{7}^\la(\tau)d\tau\\
&\leq2\cdot6^q\hc K^{q}(2K+3K_2)\|\xi\|\cdot |\Phi_1-\Phi_2|'.
\end{align*}
If $\hc$ is sufficiently small, then $J^\la$ is a
contraction. Therefore, for each $\la\in Y$, there exists a unique function $\Phi^\la\in
\bar{\Xa}$ such that \eqref{eqyjj} holds for any $(s,\xi)\in
Z_\bt$. From the one-to-one
correspondence between $\Xa$ and $\bar{\Xa}$, it follows that there exists a
unique function $\Phi^\la\in\Xa$ such that \eqref{eqyjj} holds for each $\la\in Y$ and any $(s,\xi)\in
Z_\bt$.
\end{proof}

\begin{lemma}\label{le123456}
For each $(s,\xi)\in Z_{\bt}$, let $u_i(t)=u^{\Phi^{\la_i},\la_i}(t,s,\xi),
i=1,2$, then there exists a $K_3>0$ such that
$$
\|u_1(t)-u_2(t)\|\leq K_3(h(t)/h(s))^a\mu(s)^\ve\|\la_1-\la_2\|\cdot\|\xi\|.
$$
\end{lemma}
\begin{proof}
By Lemma \ref{leaaa}, \ref{lezx}, and \ref{lemma.4}, we have
\begin{align*}
B_1^{\la_1,\la_2}(\tau):&=\|f(\tau,u_1(\tau),\Phi^{\la_1}(\tau,u_1(\tau)),\la_1)-
f(\tau,u_2(\tau),\Phi^{\la_2}(\tau,u_2(\tau)),\la_2)\|\\
&\leq\|f(\tau,u_1(\tau),\Phi^{\la_1}(\tau,u_1(\tau)),\la_1)-
f(\tau,u_1(\tau),\Phi^{\la_1}(\tau,u_1(\tau)),\la_2)\|\\
&\quad+\|f(\tau,u_1(\tau),\Phi^{\la_1}(\tau,u_1(\tau)),\la_2)-
f(\tau,u_2(\tau),\Phi^{\la_2}(\tau,u_2(\tau)),\la_2)\|\\
&\leq 6^{q+1}K^{q+1}\hc (h(\tau)/h(s))^{a(q+1)}\mu(s)^{\ve(q+1)}\\
&\quad\times
[|\la_1-\la_2|\cdot\|\xi\|^{q+1}+2\|\xi\|^q\|u_1-u_2\|_*
+\f{2}{3}|\Phi^{\la_1}-\Phi^{\la_2}|'\cdot\|\xi\|^{q+1}]
\end{align*}
and
\begin{align*}
\|\Phi^{\la_1}(s,\xi)-\Phi^{\la_2}(s,\xi)\|&\leq\il_{s}^\iy\|T(\tau,s)^{-1}Q(\tau)\|B^{\la_1,\la_2}_1(\tau)d\tau\\
&\leq
h'|\la_1-\la_2|\cdot\|\xi\|+2 h'\|u_1-u_2\|_*+(2/3)
h'|\Phi^{\la_1}-\Phi^{\la_2}|'\cdot\|\xi\|.
\end{align*}
where $h'=2\cdot3^{q+1}K^{2}\hc$. If $\hc$ is sufficiently small, let $H=h'/(1-(2/3)h')$, then
\begin{align*}
\|\Phi^{\la_1}(s,\xi)-\Phi^{\la_2}(s,\xi)\|\leq
H|\la_1-\la_2|\cdot\|\xi\|
+2H\|u_1-u_2\|_*.
\end{align*}
It follows from \eqref{eqbl} that
\begin{align*}
\|u_1(t)-u_2(t)\|&\leq
\il_s^t\|T(s,\tau)P(\tau)\|B^{\la_1\la_2}_1(\tau)d\tau\\
&\leq h'\left((1+(2/3)H)|\la_1-\la_2|\cdot\|\xi\|
+(2+(4/3)H)\|u^{\la_1}-u^{\la_2}\|_*\right)(h(t)/h(s))^a\mu(s)^\ve
\end{align*}
Thus,
$$
\|u_1-u_2\|_*\leq
[K_3/(2K)]|\la_1-\la_2|\cdot\|\xi\|,
$$
where $K_3=h'(1+2H/3)/(1-h'(1+2H/3)/K)$.
\end{proof}

We are now at the right position to establish Theorem
\ref{theoremlsm}.
\begin{proof}[{\bf Proof of Theorem \ref{theoremlsm}.}]
Sum up the above claims, we have the following conclusions.
\begin{itemize}
\item From Lemma~\ref{leaaa}, it follows that, for any $(s,\xi,\Phi,\la)\in
Z_\bt\times\bar{\Xa}\times Y$, there exists a unique function
$u(t)=u^{\Phi,\la}(t,s,\xi) \in \Omega_4$. By Lemma \ref{leaabb} ,
\ref{leopq} and the one-to-one correspondence between $\Xa$ and
$\bar{\Xa}$, for $s\geq0$ and $\xi\in
B_s\left((\bt(s)\cdot\mu(s))^{-\ve}/(2K)\right)$, there exists a
unique function $\Phi\in\Xa$ such that \eqref{eqzwt} holds for each $\la\in Y$. For
$(s,\xi)\in Z_{\bt\cdot\mu}(2K)$,  by \eqref{eqcxy}, one has
\begin{align*}
\|u(t)\|&\leq2K(h(t)/h(s))^a\mu(s)^{\ve}\f{1}{2K}(\bt(s)\cdot\mu(s))^{-\ve}\leq
(h(t)/h(s))^a\bt(s)^{-\ve}\leq\bt(s)^{-\ve},
\end{align*}
which implies that $(t,u(t))\in Z_\bt$, $t\geq s$.
Therefore, \eqref{eqkkk} holds and  $\W^\la$ is forward invariant with
respect to the semiflow $\Psi_\kappa^\la$ for each $\la\in Y$.

\item For any $(s,\xi_1), (s,\xi_2) \in Z_{\bt\cdot\mu}(2K)$, $\la\in Y$, and
$\kappa=t-s \ge 0$, by Lemma~\ref{lezx}, we have
\begin{align*}
&\lVert \Psi_\kal^\la(s,\xi_1,\Phi(s,\xi_1))-\Psi_\kal^\la(s,\xi_2,\Phi(s,\xi_2))\rVert \\
&= \lVert(t,u^{\Phi,\la}(t,s, \xi_1),\Phi^\la(t,u^{\Phi,\la}(t,s, \xi_1)))-
(t,u^{\Phi,\la}(t,s, \xi_2),\Phi^\la(t,u^{\Phi,\la}(t,s, \xi_2)))\rVert
\\
&\le3\|u^{\Phi,\la}(t,s, \xi_1)-u^{\Phi,\la}(t,s, \xi_2)\|
\le3K_1(h(t)/h(s))^a\mu(s)^\ve\|\xi_1-\xi_2\|.
\end{align*}

\item It follows from Lemma \ref{le123456} that
 for $(s,\xi)\in Z_{\beta\cdot\mu}(2K)$, $\la_1,\la_2\in Y$, and $\kappa=t-s \ge 0$, we have
\begin{align*}
&\lVert \Psi_\kal^{\la_1}(s,\xi,\Phi^{\la_1}(s,\xi))-\Psi_\kal^{\la_2}(s,\xi,\Phi^{\la_2}(s,\xi))\rVert \\
&=
\lVert(t,u_1(t),\Phi^{\la_1}(t,u_1(t))-
(t,u_2(t),\Phi^{\la_2}(t,u_2(t))\rVert
\\
&\le\|u_1(t)-u_2(t)\|+\|\Phi^{\la_1}(t,u_1(t))-\Phi^{\la_2}(t,u_2(t))\|\\
&\leq\|u_1(t)-u_2(t)\|+|\Phi^{\la_1}(t,u_1(t))-\Phi^{\la_1}(t,u_2(t))\|+\|\Phi^{\la_1}(t,u_2(t))-\Phi^{\la_2}(t,u_2(t))\|\\
&\leq3\|u_1(t)-u_2(t)\|+\|\Phi^{\la_1}-\Phi^{\la_1}|'\|u_2(t)\|\\
&\leq[3K_3+2KH(1+\f{K_3}{K})]\left(\f{h(t)}{h(s)}\right)^a\mu(s)^\ve\|\la_1-\la_2|\cdot\|\xi\|.
\end{align*}
\end{itemize}
The proof of Theorem \ref{theoremlsm} is complete.
\end{proof}

\begin{remark}
Theorem \ref{theoremlsm} includes and extends Theorem 1 in \cite{Barreira2009f}. In particular, if the parameter $\lambda$ is absent from \eqref{eqnonlinearaaa}, then Theorem \ref{theoremlsm} includes and extends some existing stable manifold theorems, for example, Theorem 1 in \cite{Barreira2005d} (nonuniform
exponential dichotomy), Theorem 2.1 in \cite{Bento2013}
(nonuniform $(\mu,\nu)$-dichotomy), and Theorem 2 in
\cite{Barreira2009e} ($\rho$-nonuniform exponential dichotomy).
\end{remark}

\section*{Acknowledgement}
The authors thank Prof. Kenneth James Palmer for his carefully reading of the manuscript, suggestive discussion, and illuminating comments, and also thank Prof. Luis Barreira for leading them into the field of the nonuniform dichotomies.

\end{document}